\documentclass[number,citesort,seceqn,dvips]{arxbj}
\usepackage{upgreek}
\usepackage{graphicx}

\aid{0}
\volume{17}
\issue{1}
\pubyear{2011}
\firstpage{347}
\lastpage{394}
\doi{10.3150/10-BEJ272}

\makeatletter
\newcommand{\eq}[1]{(\ref{#1})}
\newtheorem{lemma}{Lemma}[section]

\def\sfrac#1#2{#1/#2}
\def\vfrac#1#2{(#1)/#2}
\def\afrac#1#2{#1/(#2)}

\newcommand{\cH}{\mathcal{H}}

\newtheorem{theorem}{Theorem}[section]

\newremark{remark}{Remark}[section]
\makeatother

\begin{document}
\begin{frontmatter}

\title{Simultaneous critical values for $t$-tests in very high dimensions}
\runtitle{$t$-tests in very high dimensions}

\begin{aug}
\author[a]{\fnms{Hongyuan} \snm{Cao}\thanksref{a}\ead[label=e1]{hycao@uchicago.edu}} \and
\author[b]{\fnms{Michael R.} \snm{Kosorok}\thanksref{b}\ead[label=e2]{kosorok@unc.edu}}
\runauthor{H. Cao and M.R. Kosorok}
\address[a]{Department of Health Studies, 5841 South Maryland Avenue MC
2007, University of Chicago, Chicago, IL, 60637, USA.
\printead{e1}}
\address[b]{Department of Biostatistics and Department of Statistics
and Operations Research, 3101 Mcgavran-Greenberg Hall, CB 7420,
University of North Carolina at Chapel Hill, Chapel Hill, NC, 27599, USA.
\printead{e2}}
\end{aug}

\received{\smonth{6} \syear{2009}}
\revised{\smonth{1} \syear{2010}}

%
\begin{abstract}
This article considers the problem of multiple hypothesis testing using
$t$-tests. The observed data are assumed to be independently generated
conditional on an underlying and unknown two-state hidden model. We
propose an asymptotically valid data-driven procedure
to find critical values for rejection regions controlling the
$k$-familywise error rate ($k$-FWER), false discovery rate (FDR) and
the tail probability of false discovery proportion (FDTP) by using
one-sample and two-sample $t$-statistics. We only require a finite fourth
moment plus some very general conditions on the mean and variance of
the population by virtue of the moderate
deviations properties of $t$-statistics. A new consistent estimator for the
proportion of alternative hypotheses is developed. Simulation studies
support our theoretical results and
demonstrate that the power of a multiple testing procedure can be substantially
improved by using critical values directly, as opposed to the
conventional $p$-value approach. Our method is applied in an analysis of
the microarray data from a leukemia cancer study that involves testing
a large number of hypotheses simultaneously.
\end{abstract}

%

\end{frontmatter}

\section{Introduction}
Among the many challenges raised by the analysis of large data sets is
the problem of multiple testing. Examples include functional magnetic
resonance imaging, source detection in astronomy and microarray
analysis in genetics and molecular biology. It is now common practice
to simultaneously measure thousands of variables or features in a
variety of biological studies. Many of these high-dimensional
biological studies are aimed at identifying features
showing a biological signal of interest, usually through the
application of large-scale significance testing. The possible outcomes
are summarized
in Table~\ref{Table1}.

%
\begin{table}
\tablewidth=6.5cm
\caption{Outcomes when testing $m$ hypotheses}\label{Table1}
\begin{tabular*}{6.5cm}{@{\extracolsep{4in minus 4in}}l lll@{}}
\hline
Hypothesis & Accept & Reject & Total\\
\hline
Null true & $U$ & $V$ & $m_0$\\
Alternative true & $F$ & $S$ & $m_1$\\
Total & $W$ & $R$ & $m$\\
\hline
\end{tabular*}
\end{table}
Traditional methods that provide strong control of the familywise error
rate ($\mbox{FWER} = P(V \ge1)$) often have low power and can be
unduly conservative in many applications. One way around this is to
increase the number $k$ of false rejections one is willing to tolerate.
This results in a relaxed version of FWER, $k$-$\mbox{FWER} = P(V \ge k)$.

Benjamini and Hochberg \cite{r1} (hereafter referred to as ``BH'')
pioneered an alternative.
Define the false discovery proportion (FDP) to be the number of false
rejections divided by the number of rejections
($\mbox{FDP} = V/(R\cup 1)$). The only effect of the $R \cup1$ in the denominator is that the
ratio $V/R$ is set to zero when $R = 0$. Without loss of generality, we
treat $\mbox{FDP} = V/R$ and define the false discovery tail
probability $\mathit{FDTP}= P(V \ge\alpha R)$, where $\alpha$ is
pre-specified, based on the application. Several papers have developed
procedures for \mbox{FDTP} control. We shall not attempt a complete review
here, but mention the following:
van der Laan, Dudoit and Pollard \cite{r36} proposed an
augmentation-based procedure,
Lehmann and Romano \cite{r24} derived a step-down procedure and
Genoves and Wasserman \cite{r16} suggested an inversion-based
procedure, which is
equivalent to the procedure of \cite{r36} under mild conditions \cite{r16}.

The false discovery rate (FDR) is the expected FDP. BH provided a
distribution-free, finite-sample method for choosing a $p$-value
threshold that guarantees that the FDR is less than a target level
$\gamma$. Since this publication, there has been a considerable amount
of research on both the theory and application of FDR control.
Benjamini and Hochberg \cite{r2}
and
Benjamini and Yekutieli \cite{r3} extended the BH\ method to a class of
dependent tests.
A~Bayesian mixture model approach to obtain multiple testing procedures
controlling the FDR is considered in \cite{r14,r30,r31,r32,r33}.
Wu \cite{r39} considered the conditional dependence model under the
assumption of Donsker properties of the indicator function of the true
state for each hypothesis and derived asymptotic properties of false
discovery proportions and numbers of rejected hypotheses. A systematic
study of multiple testing procedures is given in the book \cite{r12}.
Other related work can be found in \cite{r9,r10}.

One challenge in multiple hypothesis testing is that many procedures
depend on the proportion of null hypotheses, which is not known in
reality. Estimating this proportion has long been known as a difficult
problem. There have been some interesting developments recently, for
example, the approach of \cite{r26} (see also \cite{r14,r16,r25,r23}).
Roughly speaking, these approaches are only successful under a
condition which \cite{r16} calls the ``purity'' condition.
Unfortunately, the purity condition depends on $p$-values and is hard to
check in practice.


The general framework for $k$-FWER, \mbox{FDTP}, FDR control and the
estimation of the proportion of alternative hypotheses is based on
$p$-values which are assumed to be known in advance or can be accurately
approximated. However, the assumption that $p$-values are always
available is not realistic. In some special settings, approximate
$p$-values have been shown to be asymptotically equivalent to exact
$p$-values for controlling FDR \cite{r15,r22}. However, these
approximations are only helpful in certain simultaneous error control
settings and are not universally applicable. Moreover, if the $p$-values
are not reliable, any procedures derived later are problematic.

This motivates us to propose a method to find critical values directly
for rejection regions to control $k$-FWER, \mbox{FDTP} and FDR by using
one-sample and two-sample $t$-statistics. The advantage of using $t$-tests
is that they require minimum conditions on the population, only
existence of the fourth moment, which is relatively easily satisfied by
most statistical distributions, rather than other stringent conditions
such as the existence of the moment generating
function. In addition, we approximate tail probabilities of both null
and alternative hypotheses accurately, rather than $p$-value approaches
that only consider the case under null hypotheses. Thus, a better
ranking of hypotheses is obtained. Furthermore, we propose a consistent
estimate of the proportion of alternative hypotheses which only depends
on test statistics. As long as the asymptotic distribution of the test
statistic is known under the null hypothesis, we can apply our method
to estimate this proportion, resulting in more precise cut-offs.

The BH\ procedure
controls the FDR conservatively at $\pi_0 \gamma$, where $\pi_0$ is the
proportion of null hypotheses and $\gamma$ is the targeted significance
level. If $\pi_0$ is much smaller than $1$, then the statistical power
is greatly compromised. The power we use in this paper is $\mbox{NDR} =
E[S]/m_1$, as defined in \cite{r40}. In the situation that $t$-statistics
can be used, our procedure gives a better approximation and more
accurate critical values can be obtained by plugging in the estimate of
$\pi_0$. The validity of our approach is guaranteed by empirical
process methods and recent theoretical advances on self-normalized
moderate deviations, in combination with Berry--Esseen-type bounds for
central and non-central $t$-statistics.

\begin{figure}

\includegraphics{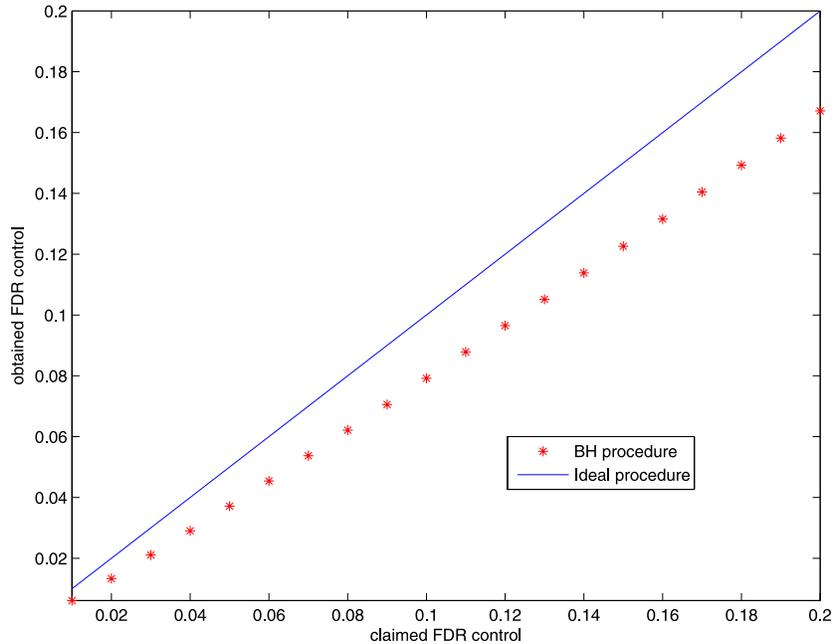}

\caption{Claimed and obtained FDR control using the BH
procedure.}\label{Figure1}
\end{figure}

To illustrate, we simulate a Markov chain, as in \cite{r34}, 
of Bernoulli variables $(H_i), i = 1, \ldots, 5000$, to
indicate the true state of each hypothesis test ($H_i = 1$ if the
alternative is true; $H_i = 0$ if the null is true). Conditional on
the indicator, observations $x_{ij}, i = 1, \ldots, 5000, j = 1,
\ldots, 80$, are generated according to the model $x_{ij} = \mu_i +
\epsilon_{ij}$. The one-sample $t$-statistic is used to perform
simultaneous hypothesis testing. Figure~\ref{Figure1} shows the plot of
10\,000 MCMC
results of the realized and nominal FDR control based on the BH\ method
for different control levels. From this plot, we can see that as the
control level increases, the BH\ procedure
becomes more and more conservative. For instance, the FDR
actually obtained is $0.167$ when the nominal level is set at $0.2$,
reflecting a significant loss in
power.
%
%

The three methods of multiple testing control we utilize are $k$-FWER,
\mbox{FDTP} and
FDR. The criterion for using $k$-FWER is, asymptotically,
%
%
\begin{equation}
P(V \ge k) \le\gamma. \label{1.3}
\end{equation}
Since we only apply our method when there are discoveries ($R >0$),
we need the \mbox{FDTP}, with a given proportion $0 < \alpha< 1$
and significance level $0 < \gamma< 1$, to satisfy, asymptotically,
%
%
\begin{equation}
P(V \ge\alpha R) \le\gamma. \label{1.1}
\end{equation}
Similarly, the criterion for using FDR is, asymptotically,
%
%
\begin{equation}
\mathit{FDR} \le\gamma\quad\mbox{or}\quad
\int_{0}^{1}P(V \ge\alpha R)\,\mathrm{d}\alpha\le\gamma. \label{1.2}
\end{equation}

The main contributions of this paper are as follows: (1)~Moderate
deviation results which
only require the finiteness of fourth moment, from which the statistic
is computed in probability
theory, are applied in multiple testing. Thus, the applicability of
this procedure is dramatically expanded: it can deal with non-normal
populations and even highly skewed populations.
(2)~The critical values
for rejection regions are computed directly, which circumvents the
intermediate $p$-value step.
(3)~An asymptotically consistent estimation
of the proportion of alternative hypotheses is developed for multiple
testing procedures under very general conditions.

The remainder of the paper is organized as follows.
In Section~\ref{sec2}, we
present the basic data structure, our goals, the procedures and theoretical
results for the one-sample $t$-test. Two-sample $t$-test results are
discussed in Section~\ref{sec3}. Section~\ref{sec4} is devoted to numerical
investigations using simulation and Section~\ref{sec5} applies our
procedure to detect significantly expressed genes in a microarray
study of leukemia cancer. Some concluding remarks and a discussion are
given in Section~~\ref{sec6}. Proofs of results from Sections~\ref
{sec2} and~\ref{sec3} are given in the \hyperref[sec7]{Appendix}.

\section{One-sample $t$-test}\label{sec2}

In this section, we first introduce the basic framework for
simultaneous hypothesis testing, followed by
our main results. Estimation of the unknown proportion of alternative
hypotheses $\pi_1$ is presented next. We conclude the section by
presenting theoretical results for the special case of completely independent
observations. This special setting is the basis for the more general
main results and is also of independent interest since fairly precise
rates of convergence can be obtained.

\subsection{Basic framework}\label{sec21}

As a specific application of multiple
hypothesis testing
in very high dimensions, we use gene expression microarray data. At the
level of single genes, researchers seek to establish
whether each gene in isolation behaves differently in a control
versus a treatment situation. If the transcripts are pairwise under
two conditions, then we can use a one-sample $t$-statistic to test for
differential
expression.

The mathematical model is
%
%
\begin{equation}\label{2.1}
X_{ij} = \mu_i + \epsilon_{ij},\qquad
1 \le j \le n, 1 \le i \le m.
\end{equation}
It should be noted that the following discussion is under this model
and does not hold in general. Here, $X_{ij}$ represents the expression
level in the $i$th gene and $j$th array.
Since the subjects are independent, for each $i$, $\epsilon_{i1},
\epsilon_{i2}, \ldots,\epsilon_{in}$ are independent random variables
with mean zero and
variance $\sigma_i^2$. The null hypothesis is $\mu_i = 0$ and the
alternative hypothesis
is $\mu_i \ne0$. For the relationship between different genes, we
propose the conditional
independence model, as follows. Let $(H_i)$ be a $\{0,1\}$-valued
stationary process and, given $(H_i)_{i=1}^{m},$ $X_{ij}, i=1, \ldots,
m$, are independently generated. The dependence is imposed on the
hypothesis $(H_i)$,
where $H_i=0$ if the null hypothesis is true and $H_i = 1$ if the
alternative is true. From
Table~\ref{Table1}, we can see that $\sum_{i = 1}^{m}H_i = m_1$ and
$\sum_{i = 1}^{m}(1 - H_i) = m_0$.
It is assumed that $(H_i)_{i=1}^{m}$ satisfy a strong law of large numbers:
%
%
\begin{equation}
{ 1 \over m} \sum_{i=1}^m H_i \to\pi_1 \in(0,1) \qquad\mbox{a.s.}
\label{dth}
\end{equation}
This condition is satisfied in a variety of scenarios, for example, the
independent case, Markov models and stationary models. Consider the
one-sample $t$-statistic
\[
T_i = \sqrt{n}\bar{X}_i/S_i,
\]
where
\[
\bar{X}_i =
\frac{1}{n}\sum_{j = 1}^{n}X_{ij},\qquad
S_i^2 =\frac{1}{n-1}\sum_{j = 1}^{n}(X_{ij} - \bar{X}_i)^2.
\]

If we use $t$ as
a cut-off, then the number of rejected hypotheses and the number of
false discoveries are, respectively,
%
%
\begin{equation}\label{2.3}
R = \sum_{i = 1}^{m} 1_{\{|T_i| \ge t\}},\qquad
V = \sum_{i = 1}^{m}(1 - H_i)1_{\{|T_i| \ge t\}}.
\end{equation}

Under the null hypothesis, it is well known that $T_i$ follows a
Student $t$-distribution with $n-1$ degrees of freedom if the sample is
from a normal distribution. Asymptotic convergence to a standard normal
distribution holds when the
population is completely unknown, provided that it has a finite fourth
moment under the null hypothesis. Moreover, under the alternative
hypothesis, $T_i$ can also
be approximated by a normal distribution, but with a shift in
location. We will show that 
%
%
\begin{eqnarray}\label{f-0}
F_0(t) &:=& P(|T_i| \ge t|H_i = 0) = P(|Z| \ge t)\bigl(1 +\mathrm{o}(1)\bigr) = 2 \bar{\Phi}(t)\bigl(1 + \mathrm{o}(1)\bigr),
\\ \label{f-1}
F_1(t) &:=& P(|T_i|
\ge t|H_i = 1) = \mbox{E}\bigl[P\bigl(\big|Z + \sqrt{n}\mu_i/\sigma_i\big| \ge t | \mu_i, \sigma_i\bigr)\bigr]\bigl(1 + \mathrm{o}(1)\bigr),
\end{eqnarray}
%
uniformly for $t=\mathrm{o}(n^{1/6})$ under some regularity conditions, where
$Z$ denotes the standard normal random variable, $\bar{\Phi}$ is the
tail probability of the standard normal distribution and the
critical values $t_{n,m}$ that control the FDTP and FDR asymptotically at
prescribed level $\gamma$ are bounded. These assumptions are fairly
realistic in practice. We do not require the critical
value for $k$-FWER to be bounded. Although we do not typically know $m_1$,
$F_0(t)$ or $F_1(t)$ in practice, we need the following theorem -- the
proof of which is given in the \hyperref[sec7]{Appendix} -- as the first step.
We will shortly extend this result, in Theorem~\ref{th3.2} below, to
permit estimation of the unknown quantities.
\begin{theorem} \label{th3.1}
Assume that $E(\epsilon_{ij}|\mu_i, \sigma_i^2) = 0$,
$\operatorname{Var}(\epsilon_{ij}|\mu_i, \sigma_i^2) = \sigma_i^2$,
$\limsup E\epsilon_{ij}^4 < \infty$,
$0 < \pi_1 < 1 - \alpha$ and (\ref{dth}) is satisfied. 
Also, assume that 
there exist $\epsilon_0 > 0$ and $c_0 > 0$ such that
%
%
\begin{equation}
P\bigl(\big|\sqrt{n}\mu_i/\sigma_i\big| \ge\epsilon_0|H_i = 1\bigr) \ge c_0 \qquad\forall
n\geq1.
\label{2.7}
\end{equation}
Let
%
%
\begin{equation}
\mu_{m}(t) = \alpha m_1 F_1(t) - (1 - \alpha) m_0 F_0(t) \label{mu}
\end{equation}
and
%
%
\begin{equation} \label{sigma}
\sigma^2_m(t) = \alpha^2 m_1 F_1(t)\bigl(1 - F_1(t)\bigr) + (1 - \alpha)^2
m_0 F_0(t) \bigl(1 - F_0(t)\bigr).
\end{equation}

\begin{enumerate}[(iii)]
\item[(i)]
If $t^{\mathit{fdtp}}_{n,m}$ is chosen such that
%
%
\begin{equation}\label{tnm}
t^{\mathit{fdtp}}_{n,m} = \inf\{t\dvtx  \mu_{m}(t) /\sigma_m(t) \ge z_\gamma\},
\end{equation}
where $z_\gamma$ is the $\gamma$th quintile of the standard normal
distribution, then
%
%
\begin{equation}
\lim_{m \to\infty} P(\mbox{FDP}\ge\alpha) = \lim_{m \to\infty}P(V
\ge
\alpha R) \le\gamma\label{t3.1-1}
\end{equation}
holds.

\item[(ii)] If $t^{\mathit{fdr}}_{n,m}$ is chosen such that
%
%
\begin{equation}\label{t3.1-2}
t^{\mathit{fdr}}_{n,m} = \inf\biggl\{t\dvtx  \frac{m_0 F_0(t)}{m_0 F_0(t) + m_1 F_1(t)}
\le\gamma\biggr\},
\end{equation}
then
%
%
\begin{equation}
\lim_{m \to\infty} \mathit{FDR} = \lim_{m \to\infty}E(V/R) \le\gamma
\label{t1.1d}
\end{equation}
holds.

\item[(iii)] If $t^{k\mbox{-}\mathit{FWER}}_{n,m}$ is chosen such that
%
%
\begin{equation}
t^{k\mbox{-}\mathit{FWER}}_{n,m} = \inf\bigl\{t\dvtx  P\bigl(\eta(t) \ge k\bigr) \le\gamma\bigr\},
\end{equation}
where $ \eta(t) \sim\operatorname{Poisson} (\theta(t)) $ and
\[
\theta(t)= m_o F_0(t),
\]
then 
%
%
\begin{equation}
\lim_{ m \to\infty} k\mbox{-FWER} = \lim_{m \to\infty} P(V \ge k)
\le\gamma\label{t3.1-3}
\end{equation}
holds.
\end{enumerate}
\end{theorem}
\begin{remark}
In the next section, we use a Gaussian approximation for $F_0(t)$ and
$F_1(t)$ for both FDTP and FDR, for which the critical values are shown
to be bounded. In this case, $m$ can be arbitrarily large, while the
critical value remains bounded. Due to sparsity, we use a Poisson
approximation for $k$-FWER, for which the critical value is no longer
bounded as $m \to\infty$, and we require $\log m = \mathrm{o}(n^{1/3})$.
\end{remark}

\subsection{Main results}

Note that in Theorem~\ref{th3.1}, there are an unknown parameter $m_1$
and unknown functions $F_0(t)$ and $F_1(t)$ involved in $\mu_m(t)$
and $\sigma_m(t)$. For practical settings, we need to estimate these
quantities.
We will begin by assuming that we have a strongly
consistent estimate of $\pi_1$ and will
then provide one such estimate in the next section.
Given $\mathcal{H}$, note that $p(t) = P(|T_i| \ge t) = (1 -
H_i)P(|T_i| \ge
t|H_i = 0) + H_i P(|T_i| \ge t |H_i = 1)$ can be estimated from the
empirical distribution $\hat{p}_m(t)$ of $\{|T_i|\}$, where
%
%
\begin{equation}
\hat{p}_m(t) = \frac{1}{m}\sum_{i = 1}^{m}I_{\{|T_i| \ge t\}},
\label{hatp}
\end{equation}
and that $P(|T_i|\geq t |H_i=0)$ is close to $P(|Z|\geq t)$ when $n$
is large, by (\ref{f-0}). The next theorem, proved in the \hyperref[sec7]{Appendix},
provides a
consistent estimate of the critical value $t_{n,m}$.
\begin{theorem} \label{th3.2}
Let
%
%
\begin{equation}\label{nv}
\nu_m(t) = \alpha\hat{p}_m(t) - 2(1 - \hat{\pi}_1) \bar{\Phi}(t)
\end{equation}
and
%
\begin{eqnarray}\label{tau-dep}
\tau^2_m(t)& = & \alpha^2 \bigl(\hat{p}_m(t) - 2(1-\hat
{\pi}_1)
\bar{\Phi}(t)\bigr) \biggl( 1- { 1 \over\hat{\pi}_1} \bigl(\hat{p}_m(t)
- 2 (1-\hat{\pi}_1)\bar{\Phi}(t)\bigr)\biggr)\nonumber\\[-8pt]\\[-8pt]
&&{} +  2(1 - \alpha)^2 (1-\hat{\pi}_1) \bar{\Phi}(t)\bigl(1-2 \bar{\Phi}(t)\bigr),
\nonumber
\end{eqnarray}
where $\hat{\pi}_1$ is a strongly consistent estimate of
$\pi_1$. Assume that the conditions of Theorem~\ref{th3.1} are
satisfied.
\begin{enumerate}[(iii)]
\item[(i)]If $\hat{t}_{n,m}^{\mathit{fdtp}}$ is chosen such that
%
%
\begin{equation}
\hat{t}_{n,m}^{\mathit{fdtp}} = \inf\biggl\{t\dvtx  \frac{\sqrt{m} \nu_m(t)}{\tau_m(t)}
\ge z_\gamma\biggr\}, \label{t3.2-1}
\end{equation}
then
%
%
\begin{equation}
|\hat{t}_{n,m}^{\mathit{fdtp}} - t_{n,m}^{\mathit{fdtp}}| = \mathrm{o}(1)\qquad\mbox{a.s.}
\end{equation}

\item[(ii)]If $\hat{t}_{n,m}^{\mathit{fdr}}$ is chosen such that
%
%
\begin{equation}
\hat{t}_{n,m}^{\mathit{fdr}} = \inf\biggl\{t\dvtx  \frac{2(1 - \hat{\pi}_1)\bar{\Phi
}(t)}{\hat{p}_m(t)} \le\gamma\biggr\}, \label{t3.2-2}
\end{equation}
then
%
%
\begin{equation}
|\hat{t}_{n,m}^{\mathit{fdr}} - t_{n,m}^{\mathit{fdr}}| = \mathrm{o}(1) \qquad\mbox{a.s.}
\end{equation}

\item[(iii)]
If $\hat{t}_{n,m}^{k\mbox{-}\mathit{FWER}}$ is chosen such that
%
%
\begin{equation}
\hat{t}_{n,m}^{k\mbox{-}\mathit{FWER}} = \inf\bigl\{t\dvtx  P\bigl(\zeta(t) \ge k\bigr)\bigr\} \le
\gamma,
\label{t3.2-3}
\end{equation}
where $ \zeta(t) \sim\operatorname{Poisson} (\bar{\theta}(t)) $ and
\[
\bar{\theta}(t) =2 m(1 - \hat{\pi}_1)\bar{\Phi}(t),
\]
then, as long as $\log m = \mathrm{o}(n^{1/3})$, we have
%
%
\begin{equation}
|\hat{t}_{n,m}^{k\mbox{-}\mathit{FWER}} - t_{n,m}^{k\mbox{-}\mathit{FWER}}| = \mathrm{o}(1) \qquad\mbox{a.s.}
\end{equation}
\end{enumerate}
\end{theorem}
\begin{remark}
This theorem deals with the general dependence case, where $(H_i)_1^m$
is assumed to follow a two-state hidden model and the data are generated
independently conditional on $(H_i)_1^m$. The proof is mainly based on
the independence case, which we present in Section~\ref{sec24} below, plus a
conditioning argument.
\end{remark}

\subsection{Estimating $\pi_1$}\label{sec23}

In the previous section, we assumed that $\hat{\pi}_1$ was a
consistent estimator of $\pi_1$. We now develop one such
estimator. By the two-group nature of multiple testing, the test
statistic is essentially a mixture of null and alternative hypotheses
with proportion as a parameter. By virtue of moderate deviations, the
distribution of $t$-statistics can be accurately approximated under both
null and alternative hypotheses. However, for the alternative
approximation, an unknown mean and variance are involved. So, we think
of a functional transformation of the $t$-statistics which has a ceiling
at $1$ to first get a conservative estimate of $\pi$ which is
consistent under certain conditions. Let $c>0$ and define $g_c(x)=
\min(|x|, c)/c$. It is easy to see that $g_c$ is a decreasing
function of~$c$, bounded by $1$,
and that the derivative ${\mathrm{d}g_c \over \mathrm{d}c}$ is bounded by $1/c$. Hence,
the function class $\{g_c\}$ indexed by $c$
is a Donsker class and thus also Glivenko--Cantelli. Let
%
%
\begin{equation}
\hat{g}_c = { 1 \over m} \sum_{i=1}^m g_c(T_i). \label{hatg}
\end{equation}
\begin{theorem} \label{p2.1}
We have
\[
\pi_1 \geq\lim_{m \to\infty, n \to\infty} \sup_{c > 0} {
\hat{g}_c - E(g_c(Z)) \over1- E(g_c(Z))} \qquad\mbox{a.s.} \label{p2.1-01}
\]
%
If, in addition, we assume that
%
%
\begin{equation}
\sqrt{n} \mu_i/\sigma_i \to\infty\qquad \mbox{for all } i
\mbox{ with } H_i= 1, i = 1, \ldots, m, \mbox{ a.s. as } n
\to\infty, \label{p2.1-02}
\end{equation}
then
\[
\pi_1 = \lim_{m \to\infty, n \to\infty} \sup_{c >0} {
\hat{g}_c - E(g_c(Z)) \over1- E(g_c(Z))} \qquad\mbox{a.s.}, \label{p2.1-03}
\]
where
\[
E(g_c(Z))= { 2 \over c \sqrt{2 \uppi}} ( 1- \mathrm{e}^{-c^2/2}) + 2\bar{\Phi}(c).
\]
%
\end{theorem}

\begin{pf}
We can write
\begin{eqnarray*} \hat{g}_c& =&
\frac{\sum_{i=1}^{m}1_{\{H_i = 0\}}}{m}\frac{\sum_{i =
1}^{m}g_c(T_i)1_{\{H_i = 0\}}}{\sum_{i=1}^{m}1_{\{H_i = 0\}}} +
\frac{\sum_{i=1}^{m}1_{\{H_i = 1\}}}{m}\frac{\sum_{i =
1}^{m}g_c(T_i)1_{\{H_i = 1\}}}{\sum_{i=1}^{m}1_{\{H_i = 1\}}}\\
&:=&\frac{m_0}{m} I
+\frac{m_1}{m} \mathit{II}.
\end{eqnarray*}

Let $\mathcal{H}= \{H_i, 1 \le i \le m\}$. Conditional on $\mathcal
{H}$, $T_i, 1
\le i \le m$, are independent random variables. We consider I first.
Let
\[
A_m(c)=\frac{\sum_{i = 1}^{m}g_c(T_i|\mathcal{H})1_{\{H_i = 0\}
}}{\sum_{i = 1}^{m}
1_{\{H_i = 0\}} } - \frac{\sum_{i = 1}^{m}E (g_c(T_i|\mathcal
{H})1_{\{H_i =
0\}}}{\sum_{i = 1}^{m} 1_{\{H_i = 0\}} },
\]
let $E$ be the infinite sequence $1_{\{H_1=0\}},1_{\{H_2=0\}},\ldots$
and let $F$ be the event that $\sum_{i=1}^m 1_{\{H_i=0\}}\rightarrow
\infty$
as $m\rightarrow\infty$. By the assumption~(\ref{dth}), we know that
$P(F)=1$. Thus,
\[
P \Bigl(\lim_{m\rightarrow\infty}\sup_{c>0}|A_m(c)|=0 \Bigr)
=E \Bigl[P \Bigl( \lim_{m\rightarrow\infty}\sup_{c>0}|A_m(c)|=0
\big|E \Bigr) \Bigr]
=1,
\]
where the second equality follows from the fact that, conditional on
$E$, the
terms in the sum are i.i.d. and thus the standard Glivenko--Cantelli
theorem applies.
Arguing similarly, based on conditioning on the sequence $1_{\{H_1=1\}
},1_{\{H_2=1\}},\ldots,$ we can also establish that
\[
\sup_{c>0} \bigg|\frac{\sum_{i = 1}^{m}g_c(T_i|\mathcal{H})1_{\{H_i = 1\}}}{\sum_{i = 1}^{m}
1_{\{H_i = 1\}} } - \frac{\sum_{i = 1}^{m}E (g_c(T_i|\mathcal{H})1_{\{H_i =
1\}}}{\sum_{i = 1}^{m} 1_{\{H_i = 1\}}} \bigg|\rightarrow0 \qquad\mbox{a.s.}
\]
Now, note that $\mathit{II} \le1$. Thus, since $m_0/m \to(1 - \pi_1)$ a.s. and
$m_1/m \to\pi_1$ a.s., we have that 
when $m \to\infty, n \to\infty,$
\begin{eqnarray*}
\hat{g}_c 
& \le& (1-\pi_1) E(g_c(Z))+ \pi_1  \qquad\mbox{a.s.} \\
& = & E(g_c(Z)) + \bigl(1 - E(g_c(Z))\bigr)\pi_1.
\end{eqnarray*}
We now have
the following lower bound for $\pi_1$:
%
%
\begin{equation}
\pi_1 \geq\lim_{m \to\infty, n \to\infty} \sup_{c > 0} {
\hat{g}_c - E(g_c(Z)) \over1- E(g_c(Z))} \qquad\mbox{a.s.} \label{pi-02}
\end{equation}

Define
\begin{eqnarray*}
\Delta_1 &:=& (1 - \pi_1)E(g_c(Z)) + \pi_1 \frac{1}{m_1}\sum_{i =
1}^{m}E(g_c(T_i) | \mathcal{H}) 1_{\{H_i = 1\}},
\\
\Delta_2 &:=& (1 - \pi_1)E(g_c(Z)) + \pi_1 \frac{\sum_{i =
1}^{m}E(g_c(Z + \sfrac{\sqrt{n}\mu_i}{\sigma_i})) 1_{\{H_i =
1\}}}{\sum_{i =1}^{m} 1_{\{H_i = 1\}}}.
\end{eqnarray*}
Letting $n \to\infty$, we have $\sup_{c>0}|\Delta_1-\Delta
_2|\rightarrow0$ a.s. Also,
\begin{eqnarray*}
\Delta_2 & = & (1 - \pi_1)E(g_c(Z))
\\
&&{}+ \pi_1 \frac{1}{\sum_{i = 1}^{m}
1_{\{H_i = 1\}}}\sum_{i =1}^{m} E\biggl(g_c\biggl(Z +
\frac{\sqrt{n}\mu_i}{\sigma_i}\biggr)\bigl(I_{\{|Z +
\sfrac{\sqrt{n}\mu_i}{\sigma_i}| \ge c\}}+ I_{\{|Z +
\sfrac{\sqrt{n}\mu_i}{\sigma_i}| < c\}}\bigr)\biggr) H_i \\
& \ge& (1 - \pi_1)E(g_c(Z)) + \pi_1 
\frac{\sum_{i = 1}^m P(|Z + \sfrac{\sqrt{n}\mu_i}{\sigma_i}| \ge
c)H_i}{\sum_{i = 1}^{m}
1_{\{H_i = 1\}}}\\
& \ge& (1 - \pi_1)E(g_c(Z)) + \pi_1 \\
& =& E(g_c(Z)) + \pi_1 \bigl( 1- E(g_c(Z))\bigr).
\end{eqnarray*}
Note that
\[
\sup_c|\hat{g}_c - \Delta_1| \to0 \qquad\mbox{a.s. as }m \to\infty, n
\to\infty.
\]
Therefore,
\[
\hat{g}_c \ge
E(g_c(Z)) + \pi_1 \bigl( 1- E(g_c(Z))\bigr) \qquad\mbox{a.s. as }m \to\infty, n \to
\infty. 
\]
Thus, we obtain
%
%
\begin{equation}
\pi_1 \le\lim_{m \to\infty, n \to\infty} \sup_{c >0}{ \hat
{g}_c - E(g_c(Z))
\over1- E(g_c(Z))}\qquad \mbox{a.s.} \label{pi-03}
\end{equation}
\end{pf}

As a consequence of this theorem, we propose the following estimate of
$\pi_1$:
%
%
\begin{equation}
\hat{\pi}_{1} := \sup_{c> 0} { \hat{g}_c - E(g_c(Z)) 
\over1- E(g_c(Z))}, \label{pihat}
\end{equation}
where
\[
E(g_c(Z))= { 2 \over c \sqrt{2 \uppi}} ( 1- \mathrm{e}^{-c^2/2}) + 2\bar{\Phi}(c).
\]

\begin{remark}
If we use $\hat{\pi}_1$, as given in (\ref{pihat}), then
Theorem~\ref{th3.2} yields a fully automated procedure to carry out multiple
hypothesis testing in very high dimensions in practical data settings.
\end{remark}

\subsection{Consistency and rate of convergence under
independence}\label{sec24}

In order to prove the main results in the general, possibly dependent,
$t$-test setting, we need results under the assumption of independence
between $t$-tests. Specifically,
we assume in this section that $(T_i, H_i), i = 1, \ldots, m$ are independent,
identically distributed random variables with $\pi_1 = P(T_i = 1)$.
This independence
assumption can also yield stronger results than the more general
setting and is of
independent interest.

The next theorem, proved in the \hyperref[sec7]{Appendix}, provides a strong consistent
estimate of the
critical value $t_{n,m}$, as well as its rate of convergence.
\begin{theorem} \label{t2.2}
Let
%
%
\begin{equation}
\nu_m(t) = \alpha\hat{p}_m(t) - 2(1 - {\pi}_1)\bar{\Phi}(t)
\label{nu-m}
\end{equation}
and
\begin{eqnarray*}\label{tau-m}
\tau_m^2(t)& = & \alpha^2\hat{p}_m(t)\bigl(1 - \hat{p}_m(t)\bigr)
+ 4 \alpha(1- {\pi}_1)\hat{p}_m(t)\bar{\Phi}(t)\nonumber\\
&&{}+ 2(1 -\pi_1)\bar{\Phi}(t)\bigl(1 - 2\alpha- 2(1 - \pi_1)\bar{\Phi}(t)\bigr).
\end{eqnarray*}
Assume the conditions of Theorem~\ref{th3.1} with (\ref{dth})
replaced by the assumption that $(T_i, H_i), i =1, \ldots, m$, are
i.i.d.~and $\pi_1 = P(T_i = 1)$. Let $\mathcal{J}=\{i\dvtx  H_i=1\}$ be
the set
that contains the indices of alternative hypotheses. Also, assume
that $\mu_i, \sigma_i$ are i.i.d.~for $i \in\mathcal{J}$.
\begin{enumerate}[(iii)]
\item[(i)]If $\hat{t}_{n,m}^{\mathit{fdtp}}$ is chosen such that
%
%
\begin{equation}
\hat{t}_{n,m}^{\mathit{fdtp}} = \inf\biggl\{t\dvtx  \frac{\sqrt{m} \nu_m(t)}{\tau_m(t)}
\ge z_\gamma\biggr\}, \label{018}
\end{equation}
then
%
%
\begin{equation}
|\hat{t}_{n,m}^{\mathit{fdtp}} - t_{n,m}^{\mathit{fdtp}}|
= \mathrm{O}\bigl(n^{-1/2} + m^{-1/2}(\log\log m)^{1/2}\bigr) \label{02.16}
\qquad \mbox{a.s.}
\end{equation}
and
%
%
\begin{equation}
|\hat{t}_{n,m}^{\mathit{fdtp}} - t_{n,m}^{\mathit{fdtp}}|= \mathrm{O}(n^{-1/2} + m^{-1/2})
\qquad\mbox{in probability.} \label{02.17}
\end{equation}
Here, $t_{n,m}^{\mathit{fdtp}}$ is the critical value defined in (\ref{t1.1a}).
\item[(ii)]If $\hat{t}_{n,m}^{\mathit{fdr}}$ is chosen such that
%
%
\begin{equation}\label{019}
\hat{t}_{n,m}^{\mathit{fdr}} = \inf\biggl\{t\dvtx  \frac{2(1 - {\pi}_1)\bar{\Phi
}(t)}{\hat
{p}_m(t)} \le\gamma\biggr\},
\end{equation}
then
%
%
\begin{equation}
|\hat{t}_{n,m}^{\mathit{fdr}} - t_{n,m}^{\mathit{fdr}}|
= \mathrm{O}\bigl(n^{-1/2} +m^{-1/2}(\log\log m)^{1/2}\bigr)\qquad \mbox{a.s.} \label{rate-1}
\end{equation}
and
%
%
\begin{equation}
|\hat{t}_{n,m}^{\mathit{fdr}} - t_{n,m}^{\mathit{fdr}}|= \mathrm{O}(n^{-1/2} + m^{-1/2})
\qquad\mbox{in probability.} \label{rate-2}
\end{equation}
Here, $t_{n,m}^{\mathit{fdr}}$ is the critical value defined in (\ref{t1.1-2}).

\item[(iii)] If $\hat{t}_{n,m}^{k\mbox{-}\mathit{FWER}}$ is chosen such that
%
%
\begin{equation}
\hat{t}_{n,m}^{k\mbox{-}\mathit{FWER}} = \inf\bigl\{t\dvtx  P\bigl(\zeta(t) \ge k\bigr)\bigr\} \le
\gamma,
\label{024}
\end{equation}
where $ \zeta(t) \sim\operatorname{Poisson} (\bar{\theta}(t)) $ and
\[
\bar{\theta}(t) = 2m(1 - {\hat{\pi}}_1)\bar{\Phi}(t),
\]
then
%
%
\begin{equation}
|\hat{t}_{n,m}^{k\mbox{-}\mathit{FWER}} - t_{n,m}^{k\mbox{-}\mathit{FWER}}| = \mathrm{O}((\log
m)^{-1/2})\qquad \mbox{a.s.}
\end{equation}
Here $t_{n,m}^{k\mbox{-}\mathit{FWER}}$ is the critical value defined in (\ref{t1.1f}).
\end{enumerate}
\end{theorem}
\begin{remark} \label{r2.01}
If $\alpha=\gamma$ in Theorem~\ref{t2.2}, then it is not difficult
to see that $\hat{t}_{n,m}^{\mathit{fdtp}} - \hat{t}_{n,m}^{\mathit{fdr}} =
\mathrm{O}(m^{-1/2})\ a.s.$ Therefore, \eq{02.16} and \eq{02.17} remain valid
with $\hat{t}_{n,m}^{\mathit{fdtp}}$ replaced by $\hat{t}_{n,m}^{\mathit{fdr}}$. This
shows that controlling FDTP is asymptotically equivalent to controlling
FDR. This
is also true in the more general dependence case. Thus,
we will focus primarily on FDR in our numerical studies.
\end{remark}
\begin{remark}
Note that $\pi_1$ is assumed to be known in order to get a precise rate
of convergence for FDTP and FDR.
If $\hat{\pi}_1$ is estimated with rate of convergence $r_n$, then the
correct convergence rate for
the ``in probability'' result for FDR and FDTP would involve an
additional term $\mathrm{O}(r_n)$ added in (\ref{02.17}) and
(\ref{rate-2}). It is unclear what the correction would be for the
almost sure rate in (\ref{02.16}) and (\ref{rate-1}).
These corrections are beyond the scope of this paper and will not be
pursued further here. Note that the rate of $\hat{\pi}_1$ is
not needed in the main results presented in Sections~\ref{sec21}--\ref{sec23}.
\end{remark}
%

\section{Two-sample $t$-test}\label{sec3}

In this section, the results of the previous section are
extended to the two-sample $t$-test setting. The estimator of the unknown
parameter $\pi_1$ remains the same as in the one-sample case, but with $T_i$
in~(\ref{hatg}) being the two-sample, rather than one-sample, $t$-statistic.
Theoretical results for the rates of convergence under independence
are also presented, as in the previous section.

\subsection{Basic set-up and results}

When two groups, such as a
control and an experimental group, are independent, which we assume
here, a natural statistic to use is the two-sample $t$-statistic. As far
as possible, we adopt the same notation as used in the one-sample case,
and we assume that~(\ref{dth}) holds.
We observe the random variables
\begin{eqnarray*}
X_{ij} = \mu_{i} + \epsilon_{ij}, \qquad 1 \le j \le n_{1}, 1 \le i \le m,\qquad
Y_{ij} = \nu_{i} + \omega_{ij}, \qquad 1 \leq j \leq n_{2}, 1 \le i \le m,
\end{eqnarray*}
with the index $i$ denoting the $i$th gene, $j$ indicating the $j$th
array, $\mu_i$ representing the mean effect for the $i$th gene from
the first group and $\nu_{i}$ representing the mean effect for the $i$th
gene from the second group. The sampling processes for the two groups
are assumed
to be independent of each other. The sample sizes $n_{1}$ and $n_{2}$
are assumed to
be of the same order, that is, $0 < b_1 \le n_{1}/n_{2} \le b_2 <
\infty
$. We will also
assume that for each $i$, $\epsilon_{i1}, \epsilon_{i2}, \ldots,
\epsilon_{in_{1}}$ are independent
random variables with mean zero and variance $\sigma_i^2$;
$\omega_{i1}, \omega_{i2}, \ldots, \omega_{in_{2}}$ are independent
random variables with mean zero and variance $\tau_i^2$. The null
hypothesis is $\mu_{i} = \nu_{i}$, the alternative hypothesis is
$ \mu_{i} \ne\nu_{i}$ and the dependence is assumed to be generated
in the
same manner as the dependence in the one-sample setting. Consider the
two-sample $t$-statistic
\begin{eqnarray*}
T_i^* &=& \frac{\bar{X}_i - \bar{Y}_i}{\sqrt{S_{1i}^2/n_1 + S_{2i}^2/n_2}},
\end{eqnarray*}
where
\begin{eqnarray*}
\bar{X}_i &=& \frac{1}{n_1}\sum_{j = 1}^{n_1}X_{ij},\qquad
\bar{Y}_i = \frac{1}{n_2}\sum_{j = 1}^{n_2}Y_{ij}, \\
S_{1i}^2 &=& \frac{1}{n_1-1}\sum_{j = 1}^{n_1}(X_{ij} -\bar{X}_i)^2,\qquad
S_{2i}^2 = \frac{1}{n_2-1}\sum_{j = 1}^{n_2}(Y_{ij}- \bar{Y}_i)^2.
\end{eqnarray*}

Then
%
%
\begin{equation}
R = \sum_{i = 1}^{m} 1_{\{|T_i^*| \ge t\}}, \qquad
V = \sum_{i = 1}^{m}(1 - H_i)1_{\{|T_i^*| \ge t\}}.
\end{equation}

The two-sample $t$-statistic is one of the most
commonly used statistics to construct confidence intervals and carry out
hypothesis testing for the difference between two means. There are
several premises underlying the use of two-sample $t$-tests. It is
assumed that the data have been derived from populations with normal
distributions. Based on the fact that $S_{1i} \rightarrow\sigma_i,
S_{2i} \rightarrow\tau_i$ a.s., with moderate violation of the
assumption, statisticians quite often recommend using the two-sample
$t$-test, provided the samples are not too small and the samples are of
equal or nearly equal size. When the populations are not normally
distributed, it is a consequence of the central limit theorem that
two-sample $t$-tests remain valid. A more refined confirmation of this
validity under non-normality based on moderate deviations is shown in
\cite{r7}. Furthermore, under the
alternative hypothesis, the asymptotic results still hold, but with a shift
in location similar to the one-sample case under certain conditions,
that is,
\begin{eqnarray*}
P(|T_i^* |\ge t|H_i = 0) &=& P(|Z| \ge t)\bigl(1 + \mathrm{o}(1)\bigr), \\
P(|T_i^*| \ge t|H_i = 1) &=& P\biggl(\bigg|Z + \frac{\mu_i - \nu_i}{B_{n_1,n_2}}\bigg| \ge t \biggr)\bigl(1 + \mathrm{o}(1)\bigr),
\end{eqnarray*}
uniformly in $t = \mathrm{o}(n^{1/6})$, where $B_{n_1, n_2}^2 = \sigma_i^2/ n_1
+ \tau_i^2/n_{2}$. Under the
assumption of (\ref{dth}), asymptotic critical values to control
FDTP, FDR and $k$-FWER are very similar to the one-sample $t$-test case
with the one-sample $t$-statistic $T_i$ replaced by the two-sample
$t$-statistic $T_i^*$. The following theorem, proved in the \hyperref[sec7]{Appendix}, is
analogous to Theorem~\ref{th3.1} and is a necessary first step.
\begin{theorem}\label{th2.3}
Assume that $E(\epsilon_{ij}|\mu_i$, $\sigma_i^2) = 0$,
$E(\omega_{ij}|\nu_i$, $\tau_i^2) = 0$, $\operatorname{Var}(\epsilon_{ij}|\mu_i,
\sigma_i^2) = \sigma_i^2$,
$ \operatorname{Var}(\omega_{ij}|\nu_i, \tau_i^2) = \tau_i^2$, $\lim\sup
E\epsilon_{ij}^4 < \infty$, $\lim\sup E\tau_{i,j}^4 < \infty$, $ 0<
\pi_1 < 1 - \alpha$ and that (\ref{dth}) is satisfied.
Assume that there exist $\epsilon_0$ and $c_0$ such that
%
%
\begin{equation}
P\biggl(\bigg|\frac{\mu_i - \nu_i}{B_{n_1,n_2}}\bigg| \ge\epsilon_0\big|H_i = 1\biggr) \ge
c_0 \qquad\mbox{for all } n_1,n_2. \label{031}
\end{equation}
The conclusions of Theorem~\ref{th3.1} then hold with the one-sample
$t$-statistic $T_i$ replaced by the two-sample $t$-statistic $T_i^*$.
\end{theorem}

\subsection{Main results}

The unknown parameter $m_1$ and functions $F_0(t)$ and $F_1(t)$ in
Theorem~\ref{th2.3} are estimated similarly as in the one-sample
case with the one-sample $t$-statistic replaced by its two-sample
counterpart. The following theorem, the proof of which is given
in the \hyperref[sec7]{Appendix}, gives our main results for two-sample $t$-tests.
\begin{theorem} \label{th2.4}
Assume that the conditions in Theorem~\ref{th2.3} are satisfied. Replace
the one-sample $t$-statistic $T_i$ by the two-sample $t$-statistic
$T_i^{*}$ in Theorem~\ref{th3.2}. Let $\hat{\pi}_1$ be a strong
consistent estimate of $\pi_1$, as in (\ref{pihat}), using the two-sample
$t$-statistic $T_i^{*}$.
\begin{enumerate}[(iii)]
\item[(i)]If $\hat{t}_{n,m}^{\mathit{fdtp}}$ is chosen such that
%
%
\begin{equation}
\hat{t}_{n,m}^{\mathit{fdtp}} = \inf\biggl\{t\dvtx  \frac{\sqrt{m} \nu_m(t)}{\tau_m(t)}
\ge z_\gamma\biggr\}, \label{t2.4-1}
\end{equation}
then
\begin{equation}
|\hat{t}_{n,m}^{\mathit{fdtp}} - t_{n,m}^{\mathit{fdtp}}| = \mathrm{o}(1)\qquad \mbox{a.s.}
\end{equation}

\item[(ii)]If $\hat{t}_{n,m}^{\mathit{fdr}}$ is chosen such that
%
%
\begin{equation}
\hat{t}_{n,m}^{\mathit{fdr}} = \inf\biggl\{t\dvtx  \frac{2(1 - \hat{\pi}_1)\bar{\Phi
}(t)}{\hat{p}_m(t)} \le\gamma\biggr\} \label{t2.4-2}
\end{equation}
then
\begin{equation}
|\hat{t}_{n,m}^{\mathit{fdr}} - t_{n,m}^{\mathit{fdr}}| = \mathrm{o}(1) \qquad\mbox{a.s.}
\end{equation}

\item[(iii)] If $\hat{t}_{n,m}^{k\mbox{-}\mathit{FWER}}$ is chosen such that
%
%
\begin{equation}
\hat{t}_{n,m}^{k\mbox{-}\mathit{FWER}} = \inf\bigl\{t\dvtx  P\bigl(\zeta(t) \ge k\bigr)\bigr\} \le
\gamma,
\label{t2.4-3}
\end{equation}
where $ \zeta(t) \sim\operatorname{Poisson} (\bar{\theta}(t)) $ and
\[
\bar{\theta}(t) = 2m(1 - \hat{\pi}_1)\bar{\Phi}(t),
\]
then, provided $\log m = \mathrm{o}(n^{1/3})$, we have
%
%
\begin{equation}
|\hat{t}_{n,m}^{k\mbox{-}\mathit{FWER}} - t_{n,m}^{k\mbox{-}\mathit{FWER}}| = \mathrm{o}(1)\qquad \mbox{a.s.}
\end{equation}
\end{enumerate}
\end{theorem}
\begin{remark}
$\hat{\pi}_1$ can be estimated via (\ref{pihat}) by using two-sample
$t$-statistics. Theorem~\ref{p2.1} is applicable in the two-sample
setting, as well as in the
one-sample case, and consistency follows. Thus, Theorem~\ref{th2.4}
gives a fully automated procedure to conduct multiple hypothesis
testing using two-sample $t$-statistics after we plug in the $\hat{\pi
}_1$ given in (\ref{pihat}).
\end{remark}

\subsection{Consistency and rate of convergence under independence}

Results for the independence setting are needed for the proofs of the
main results, as was the case for one-sample $t$-tests. We can, once
again, obtain more precise estimation compared with the general
dependence case. The following theorem, proved in the
\hyperref[sec7]{Appendix}, gives us conditions and conclusions using
two-sample $t$-statistics for controlling \mbox{FDTP} and FDR 
asymptotically, as well as rates of convergence under the assumption
that $(T_i, H_i)$ are
independent of each other for $1 \le i \le m$. Assume that $\pi_1$ is the
proportion of the alternative hypotheses among $m$ hypothesis tests,
that is, $\pi_1 = P(H_i = 1)$. Let $\mathcal{J}= \{i\dvtx  H_i = 1\}$.
\begin{theorem} \label{in-2}
Assume the conditions of Theorem~\ref{th2.3} are satisfied. Rather
than (\ref{dth}), we assume that $(T_i, H_i)$ are independent and
identically distributed. In addition, $\pi_1 = P(T_1 = 1)$ and
$\mu_i, \sigma_i$ are i.i.d.~for $i \in\mathcal{J}$. Let
%
%
\begin{eqnarray}
p(t) &=& P(|T_1^*| \ge t ),
\\
a_1(t) &=& \alpha p(t) - (1-\pi_1) P(|T_1^*|\ge t|H_1=0),
\\
b_1^2(t) & = & \alpha^2 p(t)\bigl(1 - p(t)\bigr)
+ 2 \alpha(1 - \pi_1) p(t)P(|T_1^*| \ge t|H_1 = 0)\nonumber\\
&&{}+ (1 -\pi_1)P(|T_1^*| \ge t|H_1 = 0) \bigl(1 - 2 \alpha- (1 - \pi_1)
P(|T_1^*| \ge t|H_1 = 0)\bigr),\nonumber
\\
\hat{p}_m(t) &=& \frac{1}{m}\sum_{i = 1}^{m}I_{\{|T_i^*| \ge t\}},
\\
\nu_m(t) &=& \alpha\hat{p}_m(t) - 2(1 - \pi_1) \bar{\Phi}(t),
\end{eqnarray}
and
\begin{eqnarray*}
\tau_m^2(t)& = & \alpha^2\hat{p}_m(t)\bigl(1 - \hat{p}_m(t)\bigr)
+ 4 \alpha(1- \pi_1)\hat{p}_m(t)\bar{\Phi}(t)\nonumber\\
&&{}+ 2(1 -\pi_1)\bar{\Phi}(t)\bigl(1 - 2\alpha- 2(1 - \pi_1)\bar{\Phi}(t)\bigr).
\label{tau-t2}
\end{eqnarray*}
The conclusions of Theorem~\ref{t2.2} then hold with the one-sample
$t$-statistics $T_i$ replaced by the two-sample $t$-statistics $T_i^*$.
\end{theorem}
\begin{remark}
In the above sections, we developed our theorems based on two-sided tests.
The results for the case of one-sided tests are very similar, but with
the rejection region
$\{T_i \ge t\}$ for each test. We omit the details.
\end{remark}

\section{Numerical studies}\label{sec4}

In this section, we present numerical studies based on simulated data
and compare the power of our approach with \cite{r1} (BH) and \cite{r32} (ST)
approaches using one-sample $t$-statistics. The results for
using two-sample $t$-statistics are very similar and so we omit the
details here.

\subsection{Simulation study 1}

We investigate the results for the i.i.d.~case first. Recall the model
\[
X_{ij} = \mu_i + \epsilon_{ij},\qquad
1 \le i \le m, 1 \le j \le
n.
\]
We set the signal using $\mu_i \sim \mathit{Unif}(0.5, 1)$ or $\mu_i \sim
\mathit{Unif}(-1, -0.5),$ which is of the
correct order for the standardized error term. Here, the number of
hypothesis tests is $m = 10\, 000$, which is the same for all following
simulation studies, unless otherwise noted. The proportion of
alternatives $\pi_1= 0.2$ and the error term $t(4)$ are used just to
illustrate the
asymptotic results. We vary the number of arrays $n$ from 20 to 50
to 300 to evaluate our asymptotic approximation. Empirical distributions
of \mbox{FDTP}, FDR and $k$-FWER based on $100\, 000$ repetitions are
treated as
the gold standard since they have almost negligible Monte Carlo error.
The samples are generated to evaluate our proposed method based on asymptotic
theory. Specifically, for each sample, we calculate the sample paths of
the following quantities indexed by $t$:
$\sqrt{m}\nu_m(t)/\tau_m(t)$ for studying \mbox{FDTP}, $2(1-\hat
{\pi
}_1)\bar
{\Phi}(t)/\hat{p}_m(t)$ for studying
FDR and $P(\operatorname{Poisson}(2 m(1 - \hat{\pi}_1)\bar{\Phi}(t)) \ge10)$ for studying
10-FWER (here, we choose $k = 10$ just for the purposes of
illustration). $\hat{\pi}_1$ is defined as in (\ref{pihat}).

\begin{figure}

\includegraphics{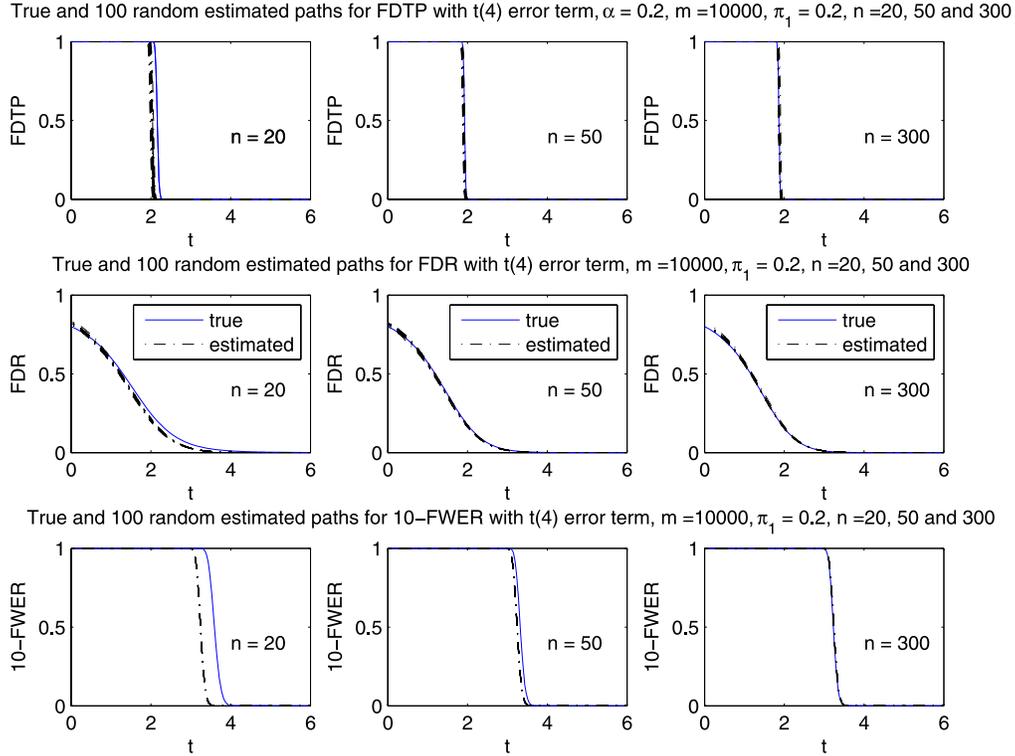}

\caption{Overlay of true and 100 random estimated sample paths with
respect to cut-off $t$ for the three procedures under differing sample
sizes.}\label{Figure1.2}
\end{figure}

Figure~\ref{Figure1.2} shows the overlay of the true path and 100
random estimated paths for \mbox{FDTP}, FDR and $k$-FWER, respectively.
As $n$
increases,
we see that the true path and
estimated paths are fairly close to each other, which, in turn,
validates our asymptotic theory.
We can see that the slopes of \mbox{FDTP} and 10-FWER are very
steep, which means a small change in the critical value results in a
large change in the level of control, while the FDR has a flatter
trend.

%

\subsection{Simulation study 2}

Under the same set-up as in the previous section, we simulate data with
different error terms: standard normal ($N(0, 1)$), Student $t$ with one
degree of
freedom (Cauchy), Student $t$ with four degrees of freedom ($t(4)$),
Student $t$ with ten degrees of freedom ($t(10)$), Laplace and
exponential. Note that, except for the Cauchy error term, all of the\vadjust{\goodbreak}
error terms satisfy the condition of finite fourth moment. Empirical
distributions of \mbox{FDTP}, FDR and $k$-FWER based on $100\, 000$ repetitions
are treated as the gold standard for obtaining true critical values.
Each scenario is repeated $1000$ times to evaluate our
proposed method for estimating the critical value based on
asymptotic theory. We control FDR at different levels (from 0.01 to
0.2) to get true
and estimated critical values.
%
%
\begin{figure}

\includegraphics{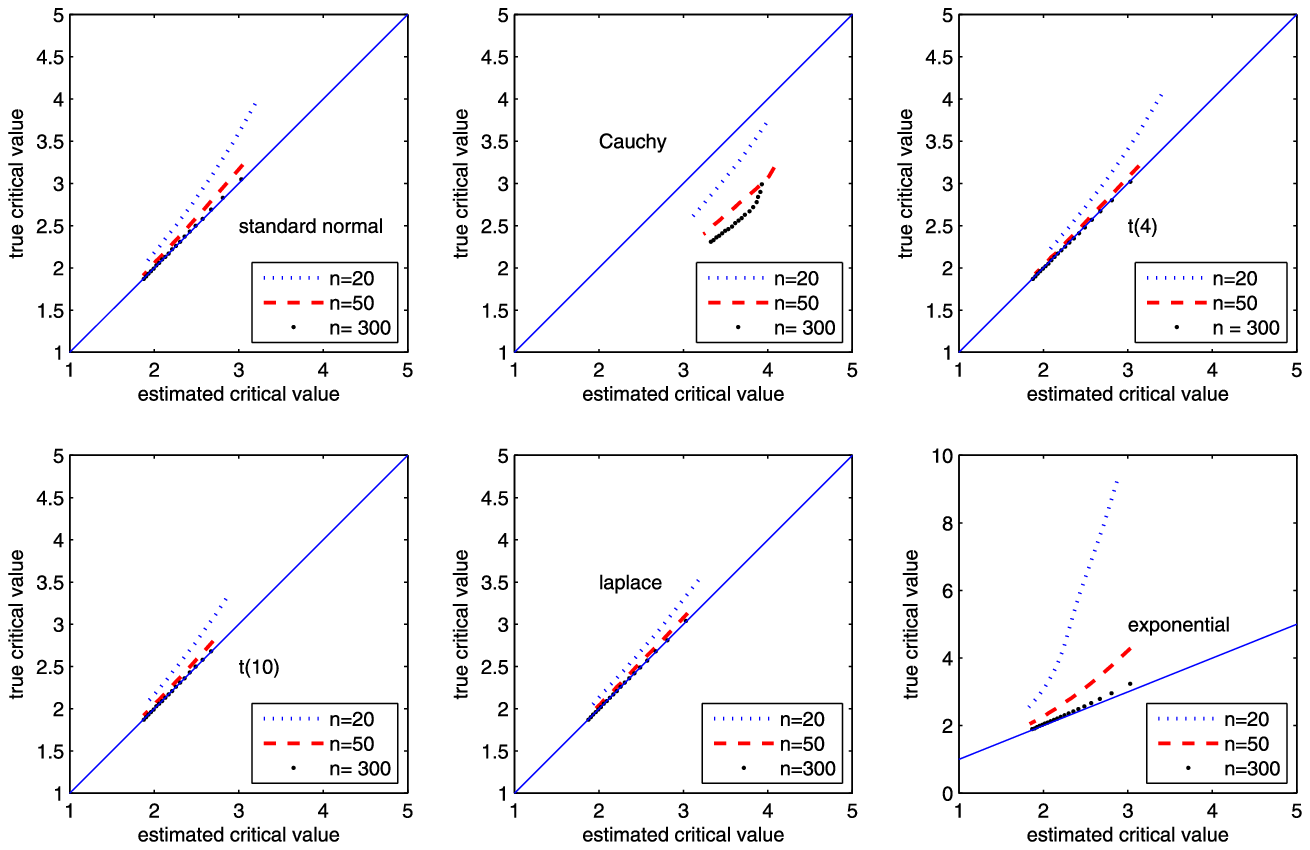}

\caption{Comparison of true and estimated critical values using FDR
for different error terms and numbers of arrays $n$.}\label{Figure2}
\end{figure}
Asymptotically, the estimated critical
value $\hat{t}$ based on our theory should be very close to the true
critical value $t$ and lie on a diagonal line of the square.
From Figure~\ref{Figure2}, the estimated
critical values $\hat{t}$ do not match the true critical value $t$ under
the Cauchy error since the Cauchy distribution does not have finite fourth
moment. For the Cauchy distribution, even the central limit theorem
does not hold since it does not have finite mean. As the number of
arrays $n$ increases, the estimated critical values $\hat{t}$ match the
true critical values
$t$ better under symmetric error terms ($N(0, 1), t(4), t(10)$ and
Laplace), but not quite so well under asymmetric errors
(e.g., exponential errors). The difficulty with the exponential error terms
suggests the value of conducting research to derive higher order approximations.
We plan to undertake this in the near future.

\subsection{Simulation study 3}

\begin{figure}

\includegraphics{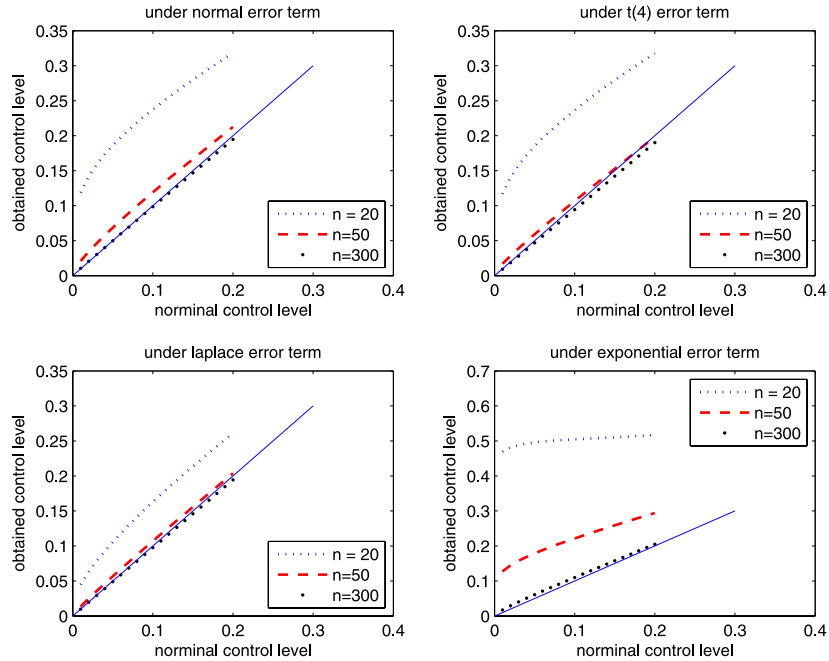}

\caption{Comparison of nominal and obtained control level for
different error terms and numbers of arrays~$n$.}\label{Figure3}
\end{figure}

The above results are from the
independent test setting.
We carried out similar simulation studies for the dependent setting and
found that
the corresponding plots are quite similar to the above results and
the same conclusions can be drawn. To see whether our proposed
method obtains the claimed level of control, we use a hidden Markov
chain to generate dependent indicators $H_i, i = 1, \ldots,
m$. Conditional on $H_i, i = 1, \ldots, m$, the data is generated
independently. The transition probability of the hidden Markov chain
is set to
\[
\pmatrix{
1 - p_1 p_1\cr
p_0 1-p_0
}
,
\]
where $p_1$ is the transition probability from 0 to 1 and $p_0$ is
the transition probability from 1 to 0. In the simulation, $p_0 = 0.8$
and $p_1 = 0.2$. Based on the limiting
stationary distribution, the alternative proportion should be $\pi_1
= p_1/(p_0 + p_1)$. Under the null hypothesis, we simulate data from
four error terms ($N(0, 1)$, $t(4)$, Laplace and exponential)
and, under the alternative hypothesis, we simulate data with mean effects
half from $\mathit{Unif}(0.1, 0,5)$ and half from $\mathit{Unif}(-0.5, -0.1),$ plus the same
four error terms. Figure~\ref{Figure3} uses FDR as the control
criterion. For
different control levels $\gamma$, we compare the claimed level of
control and the actually obtained level of control based on our
method for different numbers of arrays: small ($n = 20$), medium ($n =
50$) and large ($n = 300$).

%

%
\begin{table}[b]
\caption{Obtained control level using 10-FWER with nominal control level 0.05}\label{Table2}
\begin{tabular*}{\textwidth}{@{\extracolsep{4in minus 4in}}l l l l l@{}}
\hline
$n $ & $N(0, 1)$ & $t(4)$ & Laplace & Exponential\\
\hline
\phantom{0}20 & 0.998 (9.0e$-$05) & 0.90 (7.0e$-$03) & 0.81 (1.1e$-$02) & 1 (0) \\
\phantom{0}50 & 0.52 (1.2e$-$02) & 0.14 (9.1e$-$03) & 0.17 (1.2e$-$02) & 1 (0)\\
300           & 0.076 (3.8e$-$03) & 0.031 (2.8e$-$03) & 0.05 (2.7e$-$03) & 0.82 (4.6e$-$03) \\
\hline
\end{tabular*}
\end{table}

From Figure~\ref{Figure3}, we can see that when the number of arrays $n$
is small ($n = 20$), we do not, in general, achieve the claimed level
of control.
If we have a medium sample size ($n = 50$), the obtained level of
control is very close to the nominal level of control and the results are
almost perfect if we have a large number of arrays ($n = 300$),
even for the asymmetric exponential error term. This strongly
supports our theoretical predictions but suggests that higher order
approximations would be useful in some settings.

To see the performance of our method using 10-FWER, Table~\ref{Table2}
summarizes the control level actually obtained for different error
terms and numbers of arrays $n$ when the nominal control level is 0.05.
The obtained control level is incorrect when the number of arrays $n$ is
small, which can be deduced from the samples paths of 10-FWER given in
Figure~\ref{Figure1}. It has a very
steep slope, so when $n$ is small, the approximation is crude and
there is a noticeable difference between the estimated critical value and
the true critical value, yielding a big difference in the
control level. For large sample sizes, the obtained control level is
reasonably good because our asymptotic theory begins to take effect.
The exponential error setting appears not to perform as well as
the other error settings.

\subsection{Simulation study 4}

All previous numerical studies involve
the alternative proportion estimate $\hat{\pi}_1$ defined in (\ref
{pihat}). In this section, we investigate numerically how this estimate
is affected by number of arrays $n$ and compare with the alternative
estimate proposed by \cite{r32}. The first simulation set-up is similar
to the one in the previous section. We drew $N = 1000$ sets of data as
follows. Dependent indicators $H_i, i = 1, \ldots, m$, are generated
from a hidden Markov chain with the limiting alternative proportion
$\pi
_1 = 0.2$. Conditional on these, a vector of expected values, $\mu=
(\mu_1, \ldots, \mu_{m})$, was constructed. The expected values for the
true null hypotheses were set to $0$ with standard normal noise,
whereas the expected values for the alternative hypotheses were drawn
from $\mathit{Unif}(0.1, 0.5)$ plus standard normal noise. Correspondingly,
$1000$ replications of the proportion estimate $\hat{\pi}_1$ were
calculated using (\ref{pihat}). The root means square error (RMSE) is
given as
\[
\mathit{RMSE} = \sqrt{\frac{1}{N}\sum_{n = 1}^{N}\bigl(\hat{\pi}_1^{(n)} - \pi
_1^{(n)}\bigr)^2},
\]
where $\hat{\pi}_1^{(n)}$ is the estimate of $\pi_1$ for the $n$th
simulated data set and $\pi_1^{(n)}$ is the truth. Table~\ref{Table3}
summarizes the effect of $n$. As the number of arrays $n$ increases,
the RMSE gets smaller, which validates our asymptotic prediction.
%
%
\begin{table}[b]
\tablewidth=7cm
\caption{RMSE for $N = 1000$ estimated values of $\pi_1$}\label{Table3}
\begin{tabular*}{\tablewidth}{@{\extracolsep{4in minus 4in}}l l l l@{}}
\hline
$n $ & $20$ & $50$ & 300\\
\hline
RMSE & 0.0156 & 0.0136 & 0.0104 \\
\hline
\end{tabular*}
\end{table}

In the second simulation, we compare our proportion estimate with the
one using spline smoothing proposed by \cite{r32}. Recall the
proportion estimate $\pi_0(\lambda) = \#\{p_i > \lambda; i = 1,
\ldots,
m\}/(m(1-\lambda)).$ The smoothing approach proceeds as follows: first,
$\pi_0(\lambda)$ are calculated over a (fine) grid of $\lambda$; then,
a natural cubic spline $y$ with three degrees of freedom is fitted to
$(\lambda, \hat{\pi}_0(\lambda))$; finally, $\pi_0$ is estimated by
$\hat{\pi}_0 = y(1)$. The simulation set-up is similar to the previous
one, except that we have two groups here with $n_1 = 70$ and $n_2 =
80$. We change the alternative proportion to compare the performances
of our approach ($\pi_1^{\mathit{ck}}$) with the spline smoothing approach
($\pi
_1^{\mathit{st}}$) in Table~\ref{Table4}. They produce very similar results;
both are conservative, with less bias using our approach and less
variance using the spline smoothing approach. The advantage of our
approach is that it is computationally very fast, while the spline
smoothing approach requires that $p$-values are first obtained using
permutation, which is computationally much more intensive than our
approach (which can be computed directly from the $t$-statistics).
%
%
\begin{table}
\caption{Proportion estimate comparison}\label{Table4}
\begin{tabular*}{\textwidth}{@{\extracolsep{4in minus 4in}}l l l l l l l l l l@{}}
\hline
$\pi_1$ & $0.05$& $0.1$ & $0.15$ & $0.2$ & $0.25$ & $0.3$ & $0.35$ &$0.4$ & $0.45$\\
[6pt]
$\hat{\pi}_1^{\mathit{ck}}$ & 0.044 & 0.091 & 0.141 & 0.182 & 0.217 &0.255 & 0.289 & 0.335 & 0.365 \\
$\hat{\pi}_1^{\mathit{st}}$ & 0.041 & 0.081 & 0.125 & 0.161 & 0.195 & 0.236 &0.276 & 0.323 & 0.355\\
[6pt]
$\mathit{sd}(\hat{\pi}_1^{\mathit{ck}})$ & 0.042 & 0.043 & 0.041 & 0.040 & 0.046 &0.041 & 0.047 & 0.042 &0.038\\[1pt]
$\mathit{sd}(\hat{\pi}_1^{\mathit{st}})$& 0.039 & 0.041 & 0.036 & 0.040 & 0.041 & 0.038 &0.034 & 0.036 & 0.031\\
\hline
\end{tabular*}
\end{table}

\subsection{Comparison with BH and ST procedures}

\begin{figure}

\includegraphics{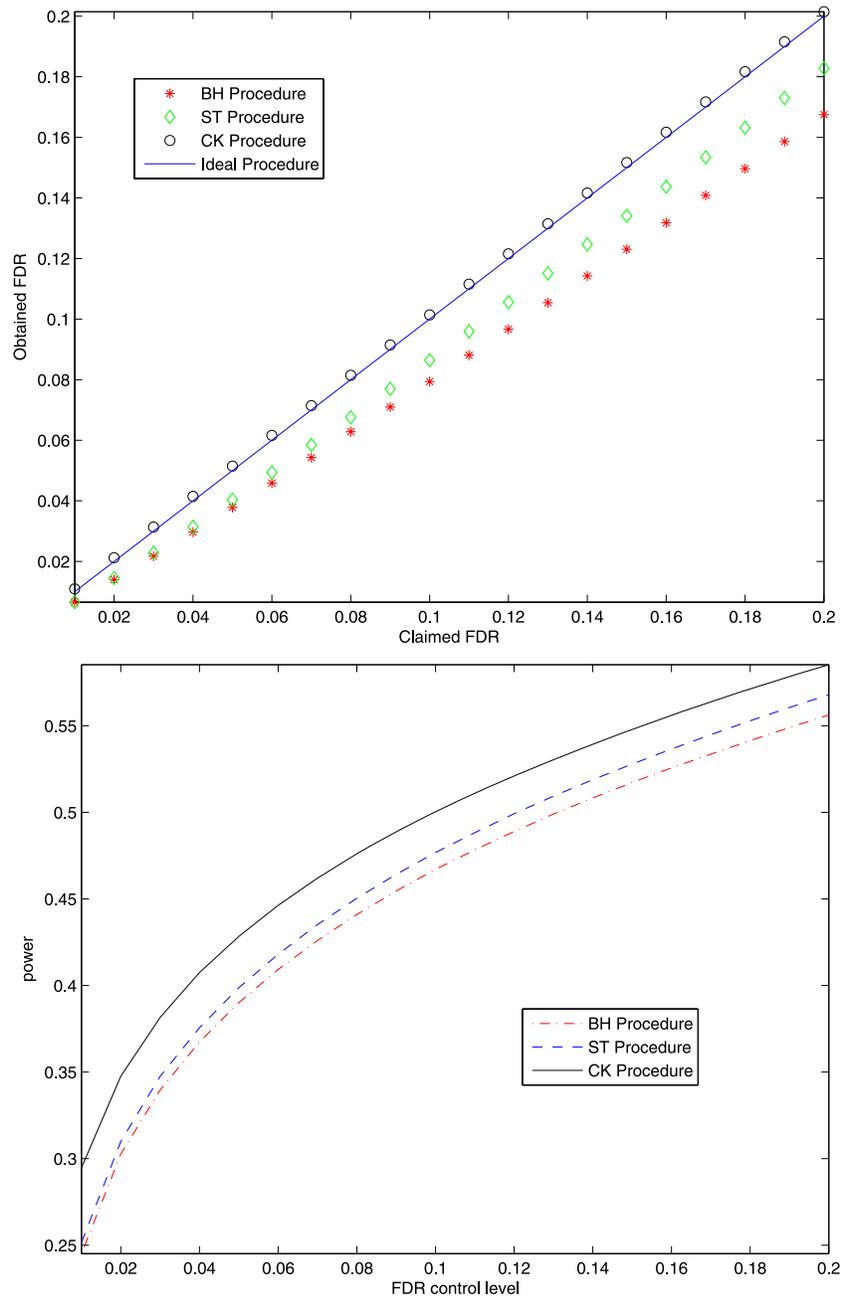}

\caption{FDR control and power comparison.}\label{Figure4}
\end{figure}

In this section, we 
compare our
approach with the BH and ST procedures under the dependence structure
described in
\cite{r39}. We also use a hidden Markov model to simulate the
indicator function $H_i, i = 1, \ldots, m$. Conditional on $H_i, i =
1, \ldots, m$, the data is generated independently. The number of
hypotheses tested $m = 5000$ and the number of arrays $n = 80$. The
data generating mechanism is otherwise the same as in the independence
case. First, we construct a one-sample $t$-statistic and apply our
procedure to obtain the critical value for the rejection region. We
then obtain $p$-values and $q$-values,
and apply the BH and ST procedures to decide which genes are significantly
expressed. We now briefly describe the BH procedure. Let $p_i$ be
the marginal \mbox{$p$-value} of the $i$th test, $1 \le i \le m$, and let
$p_{(1)} \le\cdots\le p_{(m)}$ be the order statistics of $p_1,\ldots, p_m$. Given a control level $\gamma\in(0, 1),$ let
\[
r = \max\bigl\{i \in\{0, 1, \ldots, m+1\}\dvtx  p_{(i)} \le\gamma i/m\bigr\},
\]
where $p_0 = 0$ and $p_{(m+1)} = 1.$ The BH procedure rejects all
hypotheses for which $p_{(i)} \le p_{(r)}$. If $r = 0$, then all
hypotheses are accepted. The $q$-value in \cite{r32} is
similar to the well-known \mbox{$p$-value}, except that it is a measure of
significance in terms of FDR, rather than type I error, and an estimate
of alternative proportion is plugged in, based on available $p$-values,
as described in the previous section. We revisit the motivating example
and give
a plot of the claimed FDR and actually obtained FDR by using the proposed
critical value method. 
From Figure~\ref{Figure4}, we can see that our
procedure controls the FDR at the claimed level asymptotically,
although somewhat liberally for finite samples, and
has better power at the same target FDR level compared with the BH and
ST procedures.

%

%
\begin{figure}[b]

\includegraphics{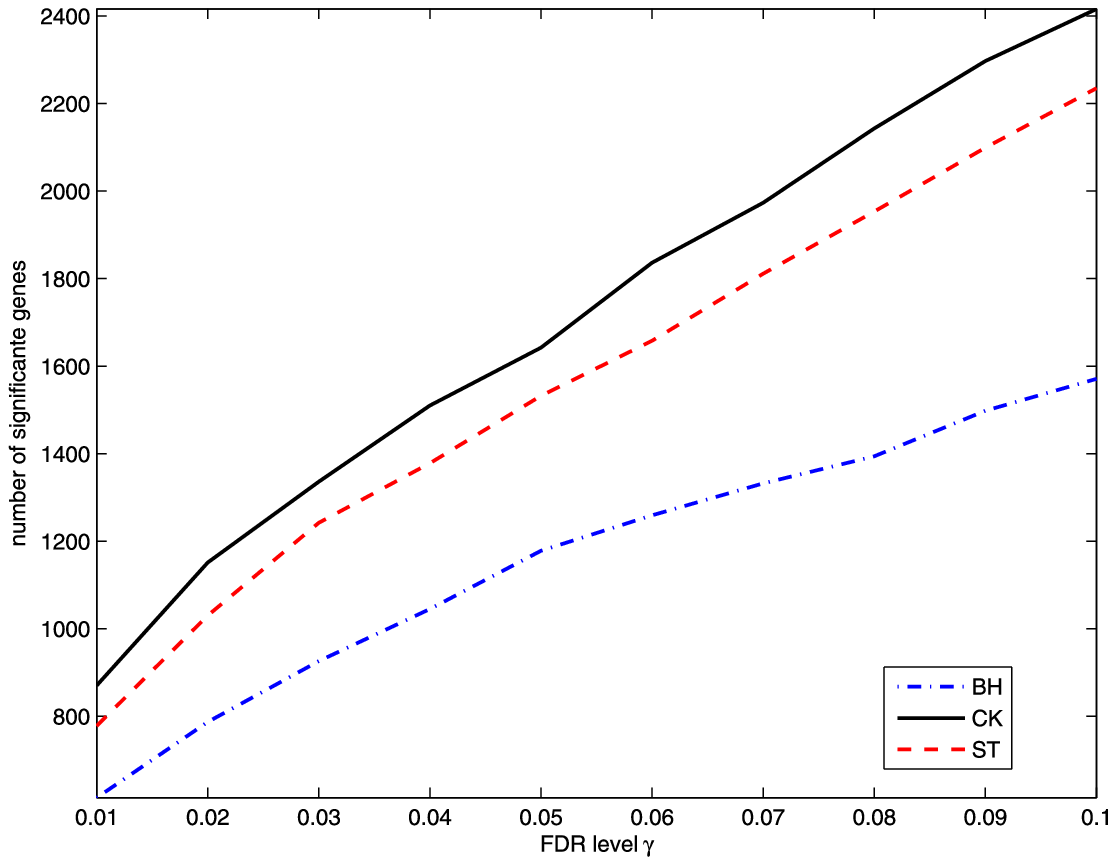}

\caption{Comparison between our (CK) procedure, the ST
procedure and the BH procedure using real data.}\label{Figure5}
\end{figure}

\section{Applications to microarray analysis}\label{sec5}

We now apply the
proposed procedure to the analysis of a leukemia cancer data set \cite
{r17} in
order to identify differentially expressed genes between AML
and ALL. For the original data, see
\url{http://www.broad.mit.edu/cgi-bin/cancer/datasets.cgi}.
In this analysis, we use the methodology developed for the dependence
case. The raw data
consist of $m = 7129$ genes and 72 samples coming from two classes:
$47$ in class
ALL (acute lymphoblastic leukemia) and 25 in class AML (acute myeloid
leukemia). Our simulation results showed
reasonable performance of the procedure for a moderate sample size in
this range. For each gene location, the two-sample $t$-statistic
comparing the $47$ ALL responses with the $25$ AML responses was
computed. Using our proposed approach for the dependent case, we
find the critical value for controlling FDR at level~$\gamma$,
\[
\hat{t}^{\mathit{fdr}}_{n,m} = \inf\biggl\{t\dvtx  \frac{2(1 - \hat{\pi}_1) \bar{\Phi}(t)}{\hat{p}_m(t)} \le\gamma\biggr\},
\]
where $\hat{p}_m = \sum_{i = 1}^{m}1_{\{|T_i| \ge t\}}/m$
and $\hat{\pi}_1$ is estimated by (\ref{pihat}). 

In Figure~\ref{Figure5}, we plot the FDR level and the number
of significantly expressed genes by our (CK) procedure, BH procedure
and the $q$-value
based Storey--Tibshirani (ST) procedure. 
From the plot, we can
see that our procedure detects the largest number of significant genes, followed
by the ST procedure and then the BH procedure, which is the most
conservative one. At FDR level $0.01$, we detected $870$ genes, the ST
procedure detected $778$ genes and the BH procedure detected $614$
genes. Using the two-sample $t$-test, similarly to the higher power of
our approach in simulation studies, we detected all of the genes that
the other two approaches detected. The BH procedure is very
conservative at the expense of power loss. The ST procedure requires
permutation to obtain $p$-values, while our procedure gets the critical
value directly and is thus faster in terms of computation. The
estimation of $\pi_1$ is $0.467$ by our procedure and $0.477$ by the ST
procedure. These results can serve as a first exploratory step for more
refined analyses concerning these significant genes. Another issue may
be that the critical value approach based on asymptotic FDR control may not
be conservative enough in some settings.

\section{Concluding remarks and discussion}\label{sec6}

We have presented a new approach for the significance analysis of
thousands of features in high-dimensional biological studies. The
approach is based on estimating the critical values of the rejection
regions for high-dimensional multiple hypothesis testing, rather than
the conventional $p$-value approaches in the literature. We developed a
detailed method that can be used to identify differentially expressed
genes in microarray experiments. The proposed procedure performs well
for large samples, reasonably well for intermediate samples and not
quite as well for small
samples, and appears to perform better than existing alternatives under
realistic sample sizes. Our method is also computationally faster than
the competing approaches. The potential for improvement in small-sample
performance motivates the need for a second-order expansion of our
theoretical work. In addition, we have proposed a new consistent
estimate of the proportion of alternative hypotheses under certain
conditions. Numerical studies demonstrate that our methodology fits the
truth well and improves the statistical power in multiple testing.
Extensions of
the current work can be pursued in several directions.

First, as stated above, the precision of the asymptotic approximations
has room for improvement in small-to-moderately-small sample sizes,
suggesting that a second-order expansion
would be valuable. Second, in the dependence case, it would be of
interest to see how the rate of convergence could be derived under
various assumptions on the form of the dependence. Thirdly,
the plug-in estimator $\pi_1$ is
consistent, but somewhat ad hoc. Complete, theoretical properties of
this estimator remain to be explored. Last, but not least, we only
considered a fixed proportion $\pi_1$ of alternative hypotheses.
It is of great interest also to consider the sparsity setting, in which
$\pi_1 \rightarrow0$ as $m \rightarrow\infty$, and to see what
patterns emerge.

\begin{appendix}
\renewcommand{\theequation}{\textup{A}.\arabic{equation}}

\section*{Appendix: Proofs of main results}\label{sec7}

Our main tools are limit theorems of empirical processes,
Berry--Esseen bounds and self-normalized moderate deviations for one-
and two-sample $t$-statistics.

\subsection{Preliminary lemmas}

We first state a non-uniform Berry--Esseen inequality for nonlinear
statistics.
\begin{lemma}[(\cite{r8})] \label{shao}
Let $\xi_1, \xi_2, \ldots, \xi_n$ be
independent random variables with $E\xi_i=0$,\break $\sum_{i=1}^n E\xi_i^2
=1$ and $E|\xi_i|^3 < \infty$. Let $W_n= \sum_{i=1}^n \xi_i$ and
$\Delta= \Delta(\xi_1, \ldots, \xi_n)$ be a measurable function of
$\{\xi_i\}$. Then
%
%
\begin{eqnarray}\label{l4.05a}
&&|P(W_n + \Delta\leq z) - \Phi(z)|\nonumber
\\
&&\quad\leq P\bigl(|\Delta|>(|z|+1)/3\bigr)
\\
&&\qquad{}+  C (|z| + 1)^{-3}
\Biggl(\|\Delta\|_2 + \sum_{i =1}^{n}\bigl(E\xi_i^2\bigr)^{1/2}\bigl(E(\Delta-\Delta_i)^2\bigr)^{1/2} + \sum_{i=1}^n
E|\xi_i|^3 \Biggr).\nonumber
\end{eqnarray}
\end{lemma}

This is \cite{r8}, Theorem 2.2, and the proof can be found there.
The next lemma provides a Berry--Esseen bound for non-central
$t$-statistics.
\begin{lemma}\label{l4.5}
Let $X, X_1, \ldots, X_n$ be i.i.d.~random variables with $E(X)=0$,
$\sigma^2= EX^2$ and $EX^4 < \infty$. Let
\[
\bar{X}= { 1 \over n} \sum_{i=1}^n X_i,\qquad
 s_n^2 = { 1 \over n-1} \sum_{i=1}^n (X_i - \bar{X})^2.
\]
Then
%
%
\begin{equation}
\bigg|P\biggl( {\sqrt{n} (\bar{X} + c) \over s_n} \leq x\biggr)
- \Phi\bigl(x - \sqrt{n} c/\sigma\bigr)\bigg| \leq K { (1+|x|) \over(1+|x-\sqrt{n} c/\sigma|)
\sqrt{n}} \label{l4.5a}
\end{equation}
for any $c$ and $x$, where $K$ is a finite constant that may depend
on $\sigma$ and $EX^4$.
\end{lemma}

\begin{pf}
Without loss of generality, assume that $x\geq0$ and $\sigma=1$.
Using
the fact that
%
%
\begin{equation}
1-|t| \leq(1+t)^{1/2} \leq1+|t| \qquad\mbox{for } t\geq-1,
\label{00}
\end{equation}
we have
%
%
\begin{equation}
x s_n = x (1+ s_n^2-1)^{1/2} \leq x (1+|s_n^2-1|) \label{l4.5-1}
\end{equation}
and
%
%
\begin{equation}
x s_n \geq x (1-|s_n^2-1|). \label{l4.5-01}
\end{equation}
Therefore,
%
%
\begin{eqnarray}\label{l4.5-2}
P\biggl( {\sqrt{n} (\bar{X} + c) \over s_n} \leq x\biggr)
& =& P\bigl( \sqrt{n} (\bar{X} + c) \leq x s_n\bigr)\nonumber\\[-8pt]\\[-8pt]
& \leq& P \bigl(\sqrt{n} \bar{X} \leq x - \sqrt{n} c + x
|s_n^2-1| \bigr).\nonumber
\end{eqnarray}

We now apply \eq{l4.05a} with $\xi_i = X_i/\sqrt{n}$, $W_n= \sqrt{n}
\bar{X}$ and
\[
z= x - \sqrt{n} c,\qquad
\Delta= -x |s_n^2-1|,\qquad
\Delta_i = -x|s_{n,i}^2 -1|,
\]
where $s_{n,i}^2$ is defined as $s_n^2$ with $0$ replacing $X_i$.

Noting that
\[
s_n^2-1 = { 1\over n-1} \Biggl(\sum_{j=1}^n (X_j^2-1) - n \bar{X}^2\Biggr) + { 1
\over n-1},
\]
\[
s_{n,i}^2-1 = { 1 \over n-1}\biggl(\sum_{j \not= i} (X_j^2-1) - n
(\bar{X}- X_i/n)^2\biggr),
\]
we have
%
%
\begin{equation}
E|s_n^2-1|^2 \leq K E X^4 / n \label{l4.5-5}
\end{equation}
and
%
%
\begin{eqnarray}\label{l4.5-6}
E(s_n^2 - s_{n,i}^2)^2 & = & { 1 \over(n-1)^2} E \bigl( (X_i^2-1) - n\bar{X}^2 +n(\bar{X}-X_i/n)^2 + 1 \bigr)^2\nonumber\\
& = & { 1 \over(n-1)^2} E \bigl( (X_i^2-1) - X_i \bigl( 2 (\bar{X}-X_i/n) +X_i/n\bigr) + 1 \bigr)^2 \nonumber\\
& \leq& { 2 \over(n-1)^2} E \bigl( 2(X_i^2-1)^2+2 + X_i^2 \bigl( 2(\bar{X}- X_i/n) + X_i/n\bigr)^2 \bigr)\\
& \leq& { 2 \over(n-2)^2} \bigl( 4 EX^4 + 6 + EX_i^2 \bigl(8 (\bar{X}-X_i/n)^2 + 2EX_i^2/n\bigr) \bigr) \nonumber\\
& \leq& K EX^4 /n^2 .\nonumber
\end{eqnarray}
It follows from \eq{l4.5-5} and \eq{l4.5-6} that
\begin{eqnarray*}
\|\Delta\|_2 &\le& K\frac{|x|\sqrt{EX^4}}{\sqrt{n}},\\
P\biggl(|\Delta| > \frac{|z|+1}{3}\biggr) &\le& K\frac{|x|\sqrt{EX^4}}{\sqrt{n}(1+ |z|)},\\
\sum_{i = 1}^{n}(E\xi_i^2)^{1/2}\bigl(E(\Delta- \Delta_i)^2\bigr)^{1/2}
&\le& K\frac{|x|\sqrt{EX^4}}{\sqrt{n}}
\end{eqnarray*}
and
\begin{eqnarray*}\sum_{i=1}^{n}E|\xi_i|^3 \le
\frac{EX^3}{\sqrt{n}}.
\end{eqnarray*}
%
Therefore, by \eq{l4.05a},
%
%
\begin{equation}\label{l4.5-9}
\big|P \bigl(\sqrt{n} \bar{X} \leq x - \sqrt{n} c + x |s_n^2-1| \bigr)
- \Phi\bigl(x -\sqrt{n} c\bigr)\big| \leq{K (1+|x|) \over(1+|x-\sqrt{n}c|)
\sqrt{n}}.
\end{equation}
Similarly,
\[
P\biggl( {\sqrt{n} (\bar{X} + c) \over s_n} \leq x\biggr)
\geq P \bigl(\sqrt{n} \bar{X} \leq x - \sqrt{n} c - x |s_n^2-1| \bigr)
\]
and
%
%
\begin{equation}\label{l4.5-10}
\big|P \bigl(\sqrt{n} \bar{X} \leq x - \sqrt{n} c - x |s_n^2-1| \bigr) -
\Phi\bigl(x -\sqrt{n} c\bigr)\big| \leq{K (1+|x|) \over(1+|x-\sqrt{n}c|)
\sqrt{n}}.
\end{equation}
This proves \eq{l4.5a}.
\end{pf}

We also need a moderate deviation for the non-central $t$-statistics, as
given in the following lemma.
\begin{lemma} \label{l4.1}
Suppose that $X, X_i, i = 1, \ldots, n,$ are independent identically
distributed random variables. Let
\begin{eqnarray*}
\bar{X} = \frac{\sum_{i = 1}^{n}X_i}{n},\qquad
s_n^2 =\frac{1}{n-1}\sum_{i = 1}^{n}(X_i - \bar{X})^2.
\end{eqnarray*}
If $X$ satisfies $E|X|^4 < \infty$, $E(X^2) = \sigma^2 > 0$ and $E(X) = 0$, then
%
%
\begin{equation}
P\biggl(\bigg|\frac{\sqrt{n} (\bar{X} + c) }{s_n} \bigg| \ge t \biggr)
= P\bigl(\big|Z+ c\sqrt{n}/\sigma\big| \ge t\bigr) \bigl(1 + \mathrm{o}(1)\bigr) \label{l4.1a}
\end{equation}
uniformly in $c$ and $t = \mathrm{o}(n^{1/6})$. Here, and in the sequel, $Z$
denotes a standard normal random variable.
\end{lemma}

\begin{pf}
When $t$ is bounded, \eq{l4.1a} follows from Lemma
\ref{l4.5}. Consider large $t$ with $t=\mathrm{o}(n^{1/6})$. We need the
following result of \cite{r37,r38}:
%
%
\begin{equation}
P\biggl(\frac{\sqrt{n} (\bar{X} + c)}{s_n} \ge t \biggr)
= \bigl(1 - \Phi\bigl(t - c\sqrt{n}/\sigma\bigr) \bigr)\bigl(1 + \mathrm{o}(1)\bigr) \label{wang}
\end{equation}
uniformly in $|c\sqrt{n}/\sigma| \le t/5$ and $t = \mathrm{o}(n^{1/6})$. We
note that following the same lines as their proof, we can see that
\eq{wang} remains valid for $- t/5\leq c \sqrt{n}/\sigma\leq t$.
We write
\[
P\biggl(\bigg|\frac{\sqrt{n} (\bar{X} + c) }{s_n} \bigg| \ge t \biggr)
= P\biggl( \frac{\sqrt{n} (\bar{X} + c) }{s_n} \ge t \biggr)
+ P\biggl( \frac{\sqrt{n}(-\bar{X} - c) }{s_n} \ge t \biggr).
\]
By \eq{wang}, the remark above and the fact that
\[
1-\Phi(t+ x) = \mathrm{o}\bigl(1-\Phi(t-x)\bigr)
\]
for $x\geq1$ (recall here that we assume $t$ is large),
\eq{l4.1a} holds for $-t \leq c\sqrt{n}/\sigma\leq t$. Now, assume
$|c|\sqrt{n}/\sigma>t$. Then, by \eq{l4.5a},
\[
\bigg|P\biggl(\bigg|\frac{\sqrt{n} (\bar{X} + c) }{s_n} \bigg| \ge t \biggr)
-P\bigl(\big|Z +c\sqrt{n}/\sigma\big|\geq t\bigr)\bigg| =\mathrm{o}(1).
\]
Since $|c|\sqrt{n} /\sigma>t$, we have $P(|Z + c\sqrt{n}/\sigma|\geq t) \geq1/2$ and hence
\[
P\biggl(\bigg|\frac{\sqrt{n} (\bar{X} + c) }{s_n} \bigg| \ge t \biggr)
=P\bigl(\big|Z +c\sqrt{n}/\sigma\big|\geq t\bigr) \bigl(1+\mathrm{o}(1)\bigr).
\]
This completes the proof of \eq{l4.1a}.
\end{pf}

The lemma below shows that $t_{n,m}$ defined in \eq{t1.1a}
under independence is bounded.
\begin{lemma}
\label{l4.2} Assume that there exist $\varepsilon_0>0$ and $c_0>0$
such that
%
%
\begin{equation}
P\bigl(\big|\sqrt{n} \mu_1/\sigma_1\big| \geq\varepsilon_0\bigr) \geq c_0. \label{l4.2a}
\end{equation}
Let $t_{n,m}$ satisfy \eq{t}. Then
%
%
\begin{equation}
t_{n,m} \leq t_0, \label{l4.2b}
\end{equation}
where $t_0$ is the solution of
%
%
\begin{equation}
\alpha\pi_1 c_0 \exp\bigl( (t_0-\varepsilon_0) \varepsilon_0\bigr) =12
(1+t_0 -\varepsilon_0).
\label{l4.2c}
\end{equation}
\end{lemma}

\begin{pf}
It suffices to show that
%
%
\begin{equation}
\sqrt{m} E\xi_1(t_0) \geq(\operatorname{var}(\xi_1(t_0)))^{1/2}
z_\gamma.
\label{l4.2-1}
\end{equation}
It is easy to see that $P(|Z + a| \ge t_0)$ is a monotone increasing
function of $a > 0$. Hence,
%
\begin{eqnarray}\label{l4.2-2}
P\bigl(\big|Z + \sqrt{n} \mu_1/\sigma_1 \big| \geq t_0\bigr)
& \geq& P\bigl(\big|Z + \sqrt{n} \mu_1/\sigma_1 \big| \geq t_0, \big|\sqrt{n}\mu_1/\sigma_1\big| \geq\varepsilon_0\bigr)\nonumber\\
& \geq& P(|Z+ \varepsilon_0| \geq t_0)P\bigl(\big|\sqrt{n} \mu_1/\sigma_1\big|\geq\varepsilon_0\bigr)\nonumber\\
& \geq& c_0 P(|Z+ \varepsilon_0| \geq t_0) \geq c_0 \bigl(1- \Phi( t_0 -\varepsilon_0)\bigr)\\
& \geq& { c_0 \over3 (1+t_0 -\varepsilon_0)} \exp\bigl(- (t_0 -\varepsilon_0)^2/2\bigr) \nonumber\\
& \geq& { c_0 \over3 (1+t_0 -\varepsilon_0)}\exp\bigl(-t_0^2 /2+(t_0-\varepsilon_0)\varepsilon_0\bigr).\nonumber
\end{eqnarray}
Here, we use the fact that
\[
\frac{1}{2}\mathrm{e}^{-x^2/2} \ge1- \Phi(x) \geq{ 1 \over\sqrt{2 \uppi}(1+x)} \mathrm{e}^{-x^2/2}
\qquad\mbox{for }x\geq0.
\]
Under the null hypothesis $H_1 = 0$, which corresponds to $\mu_i =
0$, we apply Lemma~\ref{l4.1} and obtain
%
%
\begin{equation}
P(|T_1| \ge t|H_1 = 0) = P(|Z| \geq t) \bigl(1 + \mathrm{o}(1)\bigr) \label{4.2}
\end{equation}
uniformly in $t = \mathrm{o}(n^{1/6})$.

Under the alternative hypothesis $H_1 = 1$, we apply Lemma
\ref{l4.1} to $X_{ij}- \mu_i$ and obtain
%
\begin{eqnarray}\label{4.3}
P(|T_1| \ge t|H_1=1)
&= &P\bigl( \big|\sqrt{n} ( \bar{X}_1 - \mu_1 + \mu_1)/s_1\big| \geq t|H_1=1\bigr)\nonumber\\
&= & E [P(|Z+ \sqrt{n} \mu_1/\sigma_1)| \geq t | \mu_1, \sigma_1)]\bigl(1+\mathrm{o}(1)\bigr)\\
&=& P\bigl(\big|Z + \sqrt{n} \mu_1/\sigma_1\big| \geq t\bigr) \bigl(1+\mathrm{o}(1)\bigr)\nonumber
\end{eqnarray}
uniformly in $t = \mathrm{o}(n^{1/6})$.

Also, note that
%
\begin{eqnarray}\label{4.03}
P(|T_1| \geq t)
& =& P(|T_1| \geq t,H_1=0) + P(|T_1| \geq t,H_1 =1) \nonumber\\
& =& (1-\pi_1) P(|T_1| \geq t |H_1=0) + \pi_1 P(|T_1| \geq t |H_1 =1)\nonumber\\[-8pt]\\[-8pt]
& =& (1-\pi_1)P(|Z| \geq t) \bigl(1 + \mathrm{o}(1)\bigr)\nonumber
\\
&&{}+ \pi_1 P\bigl(\big|Z + \sqrt{n}\mu_1/\sigma_1 \big| \geq t\bigr) \bigl(1+ \mathrm{o}(1)\bigr).\nonumber
\end{eqnarray}

By \eq{4.7}, \eq{4.2}, \eq{4.03} and \eq{l4.2-2},
%
\begin{eqnarray}\label{l4.2-3}
E\xi_1(t_0) &=& \alpha(1-\pi_1)P(|Z| \geq t_0) \bigl(1 + \mathrm{o}(1)\bigr)
+ \alpha\pi_1 P\bigl(\big|Z + \sqrt{n} \mu_1/\sigma_1 \big| \geq t_0\bigr)
\bigl(1+\mathrm{o}(1)\bigr)\nonumber
\\
&&{}- (1-\pi_1)P(|Z| \geq t_0) \bigl(1+\mathrm{o}(1)\bigr)\nonumber\\
& \geq& \alpha\pi_1 { c_0 \over6 (1+t_0 -\varepsilon_0)}\exp
\bigl(-t_0^2 /2+(t_0-\varepsilon_0) \varepsilon_0\bigr)- 2 P(Z \geq t_0) \nonumber\\[-8pt]\\[-8pt]
& \geq& { \alpha\pi_1 c_0 \over6 (1+t_0 -\varepsilon_0)}\exp
\bigl(-t_0^2/2 +
(t_0-\varepsilon_0) \varepsilon_0\bigr)
-\mathrm{e}^{-t_0^2/2}\nonumber \\
& =& \mathrm{e}^{-t_0^2/2} \biggl( {\alpha\pi_1 c_0
\over6 (1+t_0 -\varepsilon_0)}\exp\bigl( (t_0-\varepsilon_0) \varepsilon
_0\bigr) -1 \biggr) \nonumber\\
& =& \mathrm{e}^{-t_0^2/2},\nonumber
\end{eqnarray}
by \eq{l4.2c} and the
definition of $t_0$. It is easy to see that $E\xi_1^2 \leq1$ and
$\operatorname{var}(\xi_1(t_0)) \leq1$ in particular. Thus, by \eq{l4.2-3},
%
%
\begin{equation}
{\sqrt{m} E\xi_1(t_0) \over(\operatorname{var}(\xi_1(t)))^{1/2}}
\geq
\sqrt{m} \mathrm{e}^{-t_0^2/2} \geq z_\gamma, \label{418}
\end{equation}
provided that $m$ is large enough. This proves \eq{l4.2-1}.
\end{pf}

The following i.i.d.~results are essential for the general results.
\begin{lemma} \label{th2.1}
Assume the conditions of Theorem~\ref{th3.1} with (\ref{dth})
replaced by the assumption that $(T_i, H_i), i =1, \ldots, m$ are
i.i.d. and $\pi_1 = P(T_i = 1)$. Let $\mathcal{J}=\{i\dvtx  H_i=1\}$ be
the set
that contains the indices of alternative hypotheses. Also, assume
that $\mu_i, \sigma_i$ are i.i.d. for $i \in\mathcal{J}$. Let
%
%
\begin{eqnarray}
p(t) &=& P(|T_1| \ge t ),
\\ \label{a1}
a_1(t) &=& \alpha p(t) - (1-\pi_1) F_0(t)
\end{eqnarray}
and
%
%
\begin{eqnarray}\label{b1}
b_1^2(t) &=& \alpha^2 p(t)\bigl(1 - p(t)\bigr) + 2 \alpha(1 - \pi_1)p(t)F_0(t)\nonumber
\\[-8pt]\\[-8pt]
&&{}+ (1 -\pi_1)F_0(t) \bigl(1 - 2 \alpha- (1 - \pi_1) F_0(t)\bigr).\nonumber
\end{eqnarray}
\begin{enumerate}[(iii)]
\item[(i)]
If $t^{\mathit{fdtp}}_{n,m}$ is chosen such that
%
%
\begin{equation}
t^{\mathit{fdtp}}_{n,m} = \inf\bigl\{t\dvtx  \sqrt{m} a_1(t) /b_1(t) \ge z_\gamma\bigr\},
\label{t1.1a}
\end{equation}
then
%
%
\begin{equation}
\lim_{m \to\infty} P(\mbox{FDP}\ge\alpha) = \lim_{m \to\infty}P(V
\ge
\alpha R) \le\gamma\label{t1.1b}
\end{equation}
holds.

\item[(ii)] If $t^{\mathit{fdr}}_{n,m}$ is chosen such that
%
%
\begin{equation}
t^{\mathit{fdr}}_{n,m} = \inf\biggl\{t\dvtx  \frac{(1 - \pi_1)F_0(t)}{
p(t)} \le\gamma\biggr\}, \label{t1.1-2}
\end{equation}
then
\begin{equation}
\lim_{m \to\infty} \mathit{FDR} = \lim_{m \to\infty}E(V/R) \le\gamma
\label{t1.1d}
\end{equation}
holds.

\item[(iii)] If $t^{k\mbox{-}\mathit{FWER}}_{n,m}$ is chosen such that
%
%
\begin{equation}
t^{k\mbox{-}\mathit{FWER}}_{n,m} = \inf\bigl\{t\dvtx  P\bigl(\eta(t) \ge k\bigr) \le\gamma\bigr\},
\label{t1.1f}
\end{equation}
where $ \eta(t) \sim\operatorname{Poisson} (\theta(t)) $ and
\[
\theta(t)= m(1 - \pi_1)F_0(t),
\]
then
%
%
\begin{equation}
\lim_{ m \to\infty} k\mbox{-FWER} = \lim_{m \to\infty} P(V \ge k)
\le\gamma\label{t1.1g}
\end{equation}
holds.
\end{enumerate}
\end{lemma}

\begin{pf}
We first prove the i.i.d.~case for one-sample $t$-statistics.
By \eq{2.3},
\begin{eqnarray*}\alpha R - V & =& \alpha\sum_{i =
1}^{m}I_{\{|T_i| \ge t\}} -
\sum_{i = 1}^{m}(1 - H_i)I_{\{|T_i| \ge t\}}\nonumber\\[1pt]
& =& \sum_{i = 1}^{m}(H_i + \alpha-1) I_{\{|T_i| \ge t\}} \nonumber\\[1pt]
& =&\sum_{i=1}^m \alpha I_{\{|T_i| \ge t\}} I_{\{H_i =1\}}
+\sum_{i=1}^m (\alpha-1)I_{\{|T_i| \ge t\}} I_{\{H_i =0\}}\nonumber\\[1pt]
& =& \sum_{i=1}^m \alpha I_{\{|T_i| \ge t\}} \bigl(1- I_{\{H_i =0\}}\bigr)
+\sum_{i=1}^m (\alpha-1)I_{\{|T_i| \ge t\}} I_{\{H_i =0\}} \nonumber
\\[1pt]
& =& \sum_{i=1}^m \bigl(\alpha I_{\{|T_i|\geq t\}} - I_{\{|T_i| \ge t\}}I_{\{H_i=0\}}\bigr)\nonumber\\[1pt]
& =& \sum_{i=1}^m \xi_i,
\end{eqnarray*}
where
\[
\xi_i :=\xi_i(t) = \alpha I_{\{|T_i|\geq t\}} - I_{\{|T_i| \geq t\}}
I_{\{H_i =0\}}
\]
is obviously a Donsker class indexed by $t$ \cite{r21}.
Hence,
%
%
\begin{equation}
P(V \ge\alpha R) = P\Biggl(\sum_{i=1}^m \xi_i(t) \le0\Biggr). \label{4.5}
\end{equation}
Note that since $\xi_i$ are independent random variables, we can
apply the uniform central limit theorem to choose $t$ so that
%
%
\begin{equation}
P\Biggl(\sum_{i=1}^m \xi_i(t) \le0\Biggr) \leq\gamma. \label{4.6}
\end{equation}
To this end, we need the mean and variance of $\xi_i$.
Without loss of generality, we use $\xi_1$ as an example, since
$\xi_i$ are i.i.d.~random variables. Thus,
%
\begin{eqnarray}\label{4.7}
E \xi_1 & = & \alpha P(|T_1| \geq t) - P(|T_1| \geq t, H_1 =0)\nonumber\\[1pt]
& =& \alpha P(|T_1|\geq t) - P(H_1=0) P(|T_1| \geq t |H_1=0)\\[1pt]
& =& \alpha P(|T_1|\geq t) - (1-\pi_1) P(|T_1| \geq t |H_1=0).\nonumber
\end{eqnarray}
Similarly,
%
\begin{eqnarray} \label{4.8}
E \xi_1^2 & = & E\bigl(\alpha^2 I_{\{|T_1|\geq t\}} + (1-2\alpha)I_{\{|T_1|\geq t\}}I_{\{H_1=0\}}\bigr)\nonumber\\[-9pt]\\[-9pt]
& =& \alpha^2 P(|T_1|\geq t) + (1-2 \alpha) (1-\pi_1) P(|T_1| \geq t|H_1=0)\nonumber
\end{eqnarray}
and
%
\begin{eqnarray}\label{4.9}
\operatorname{var}(\xi_1) & = & E\xi_1^2 - (E\xi_1)^2\nonumber\\[-1pt]
& =& \alpha^2 P(|T_1|\geq t) + (1-2 \alpha) (1-\pi_1) P(|T_1| \geq t|H_1=0)\nonumber\\[-1pt]
&&{}-\{ \alpha P(|T_1|\geq t) - (1-\pi_1) P(|T_1| \geq t |H_1=0)\}^2\nonumber
\\[-8pt]\\[-8pt]
& =& \alpha^2 P(|T_1|\geq t) \bigl(1- P(|T_1| \geq t)\bigr)\nonumber\\[-1pt]
&&{}+ (1-\pi_1) P(|T_1| \geq t |H_1=0) \bigl(1 -2\alpha- (1-\pi_1) P(|T_1|\geq t|H_1=0)\bigr)\nonumber\\[-1pt]
&&{}+ 2\alpha(1-\pi_1) P(|T_1|\geq t) P(|T_1| \geq t |H_1=0).\nonumber
\end{eqnarray}
Now, define
%
%
\begin{equation}
t_{n,m} = \inf\biggl\{t\dvtx  { \sqrt{m} E\xi_1(t) \over
(\operatorname{var}(\xi_1(t)))^{1/2} } \geq z_\gamma\biggr\}. \label{t}
\end{equation}
By Lemma~\ref{l4.2}, $t_{n,m}$ is bounded and hence the uniform central
limit theorem yields
%
\begin{eqnarray}\label{4.10}
P\Biggl(\sum_{i=1}^m \xi_i(t_{n,m}) \le0\Biggr)
& =& P \biggl( {\sum_{i=1}^m (\xi_i (t_{n,m}) - E\xi_i(t_{n,m}))\over
(\sum_{i=1}^m \operatorname{var}(\xi_i(t_{n,m})))^{1/2}}\nonumber
\\[-1pt]
&\le&- {\sum_{i=1}^m
E\xi_i(t_{n,m})\over
(\sum_{i=1}^m \operatorname{var}(\xi_i(t_{n,m})))^{1/2}} \biggr)\nonumber
\\[-9pt]\\[-9pt]
& \leq& P \biggl( {\sum_{i=1}^m (\xi_i (t_{n,m}) -
E\xi_i(t_{n,m}))\over
(\sum_{i=1}^m \operatorname{var}(\xi_i(t_{n,m})))^{1/2}} \le-
z_\gamma\biggr)\nonumber\\[-1pt]
& \to& \Phi(-z_{\gamma}) = \gamma.\nonumber
\end{eqnarray}
This proves \eq{t1.1b}.

Note that
\begin{eqnarray*}
\mathit{FDR} & = & \int_0^1 P(\mathit{FDTP}\ge x)\,\mathrm{d}x \\
& = & \int_0^1 P(V \ge x R)\,\mathrm{d}x\\
& = & \int_0^1 P\Biggl(\sum_1^m \xi_i \le0\Biggr)\,\mathrm{d}x\\
& = & \int_0^1 P\Biggl(N(0,1) \le\frac{-\sqrt{m}E\xi_1}{\sqrt{\operatorname{Var}\xi_1}}\Biggr)\,\mathrm{d}x.
\end{eqnarray*}
Letting $m \rightarrow+\infty,$ $P(N(0,1) \le
-\sqrt{m}E\xi_1/\sqrt{\operatorname{Var}\xi_1})$ is either 0 or 1, depending on the
sign of $E\xi_1. $ Thus, the range of $x$ that makes this probability
1 satisfies
\[
E\xi_1 = xP(|T_1| \ge t) - (1 - \pi_1)P(|T_1| \ge t|H_1 = 0) < 0
\]
and the corresponding $x < (1 - \pi_1)P(|T_1| \ge t|H_1 = 0)/P(|T_1|
\ge
t).$ In order to control FDR at level $\gamma$, we require
\[
\frac{(1 - \pi_1)P(|T_1| \ge t|H_1 = 0)}{P(|T_1| \ge t)} \le\gamma.
\]
This proves \eq{t1.1-2}.

For the $k$-FWER, we use the characteristic function method. Let $\eta_i
= (1 - H_i)I_{\{|T_i| \ge t\}}$,
\begin{eqnarray*}
E \mathrm{e}^{\mathrm{i}s \sum_{i = 1}^{m} \eta_i}& = & \prod_{i = 1}^m E\mathrm{e}^{\mathrm{i}s \eta_i}\\
& = & \prod_{i = 1}^m [\mathrm{e}^{\mathrm{i}s}(1 - \pi_1)F_0 + 1 - (1 - \pi_1)F_0]\\
& = & \biggl[1 + \frac{1}{m} m(1 - \pi_1)F_0(\mathrm{e}^{\mathrm{i}s} - 1)\biggr]^m\\
& \to& \mathrm{e}^{\lambda(\mathrm{e}^{\mathrm{i}s} - 1)},
\end{eqnarray*}
where $m_0F_0 \to\lambda$ as $m \to\infty$ and $\lambda$ is the
parameter for the Poisson distribution such that
\[
P\bigl(\operatorname{Poiss}(\lambda) \ge k\bigr) \le\gamma.
\]
\upqed
\end{pf}

The following functional central limit theorem 
is needed in the proof of Theorem~\ref{th3.1}:
\begin{lemma} \label{lem7.8}
Suppose the triangular array $\{f_{ni}(\omega, t), i = 1, \ldots,
m_n, t
\in T\}$
consists of independent processes within rows and is almost measurable Suslin analytic set (AMS)(see page~25 in \cite{r21}). Let\vspace*{-1pt}
%
%
\begin{equation}
X_n(\omega, t) \equiv\sum_{i = 1}^{m_n}[f_{ni}(\omega, t) - E f_{ni}(\cdot, t)].
\end{equation}
Assume:
\begin{enumerate}[(A)]
\item[(A)] the $\{f_{ni}\}$ are manageable, with envelopes $\{F_{ni}\}
$ which are also
independent within rows;
\item[(B)] $H(s, t) = \lim_{n \rightarrow\infty} EX_n(s)X_n(t)$
exists for every $s, t \in T$;
\item[(C)] $\limsup_{n \rightarrow\infty} \sum_{i = 1}^{m_n}E^*
F^2_{ni} < \infty$;
\item[(D)] $\lim_{n \rightarrow\infty}\sum_{i = 1}^{m_n}E^*F_{ni}^2 1
\{F_{ni} > \epsilon\} = 0$ for each $\epsilon> 0$;
\item[(E)] $\rho(s, t) = \lim_{n \rightarrow\infty} \rho_n(s,t), $ where
\[
\rho_n(s,t) \equiv\Biggl(\sum_{i = 1}^{m_n}E|f_{ni}(\cdot, s) - f_{ni}(\cdot, t)|^2\Biggr)^{1/2}
\]
exists for every $s, t \in T$ and, for all deterministic sequences $\{
s_n\}$ and $\{t_n\}$ in $T$,
if $\rho(s_n, t_n) \rightarrow0$, then $\rho_n(s_n, t_n) \rightarrow0.$
\end{enumerate}
Then $X_n$ converges weakly on $l^{\infty}(T)$ to a tight mean-zero
Gaussian process
$X$ concentrated on $\mathit{UC}(T, \rho)$, with covariance $H(s,t).$
\end{lemma}

\begin{pf}
The definitions involved in this lemma and the proof can be found in
\cite{r21}, Theorem~11.16.
Below, we verify that, conditional on $\mathcal{H}$, $f_{ni} (\omega,
t) =
\xi
_{i}(\omega, t)/\sqrt{m}$ satisfy the
conditions in Lemma~\ref{lem7.8}. Since $\xi_i(\omega,t)$ is the
difference between
two monotone bounded functions, it is clear that,
conditional on $\mathcal{H}$, $\xi_{i}(\omega, t)/\sqrt{m}$ is $\mathit{AMS}$,
manageable and has
envelopes $\alpha/\sqrt{m}$. Also,
\begin{eqnarray*}
EX_n(s)X_n(t)& = &E E[X_n(s)X_n(t)|\mathcal{H}]\\
& = & E E \biggl[\frac{\sum_{i = 1}^{m}(\xi_i(s)|\mathcal{H}- E \xi_i(s)|\mathcal{H})}{\sqrt{m}}
\frac{\sum_{j = 1}^{m}(\xi_j(t)|\mathcal{H}- E \xi_j(t)|\mathcal{H})}{\sqrt{m}}\biggr]
\\
& = & E E \frac{\sum_{i = 1}^{m}(\xi_i(s)|\mathcal{H}- E \xi
_i(s)\mathcal{H})(\xi
_i(t)|\mathcal{H}- E \xi_i(t)\mathcal{H})}{m}\\
& = & \frac{1}{m} E \sum_{i = 1}^{m}E(\xi_i(s)|\mathcal{H})(\xi
_i(t)|\mathcal{H})
- \sum
_{i = 1}^{m}E(\xi_i(s)|\mathcal{H})E(\xi_i(t)|\mathcal{H})\\
& = & \frac{1}{m} E \sum_{i = 1}^{m}\bigl(\alpha^2 H_i + (1 - \alpha)^2(1 -H_i)\bigr)EI_{\{|T_i| \ge t\cup s | \mathcal{H}\}}\\
&&{} -  \sum_{i = 1}^{m}[\alpha H_i + (1 - \alpha)(1 - H_i)]^2E I_{\{
|T_i| \ge s \mathcal{H}\}} E I_{\{|T_i| \ge t | \mathcal{H}\}}\\
& = & \frac{1}{m} E \sum_{i = 1}^{m}\bigl(\alpha^2 H_i F_1(t \cup s) + (1 -
\alpha)^2(1 - H_i)F_0(t \cup s)\bigr) \\
&&{} -  \sum_{i = 1}^{m}[\alpha^2 H_i + (1 - \alpha)^2(1 - H_i)][H_i
F_1(s)+ (1 - H_i)F_0(s)]
\\
&&\hspace*{26pt}{}\times[H_i F_1(t) + (1 - H_i)F_0(t)]\\
& = & \frac{1}{m} E \sum_{i = 1}^{m} \bigl[\alpha^2 H_i \bigl(F_1(t\cup s) -
F_1(t)F_1(s)\bigr)
\\
&&\hspace*{35pt}{}+ (1 - \alpha)^2(1 - H_i)\bigl(F_0(t\cup s) - F_0(t)F_0(s)\bigr)\bigr]\\
& \to& \pi_1 \alpha^2 \bigl(F_1(t\cup s) - F_1(t)F_1(s)\bigr) + (1 - \pi_1)
(1 -
\alpha)^2 \bigl(F_0(t\cup s) - F_0(t)F_0(s)\bigr)\\
& \equiv& H(s, t),
\end{eqnarray*}
which is the same as $q^2(t)$ when $s = t$.
(C) is easily satisfied. For all $\epsilon> 0$, there exists an
$N_0$ such that $\alpha/ N_0 < \epsilon$, so $\lim_{m \rightarrow
\infty
}\sum_{i = 1}^{m}E \alpha^2/m 1 \{\alpha/\sqrt{m} > \epsilon\} = \lim_{m
\rightarrow\infty}\sum_{i = 1}^{N_0 - 1} \alpha^2/m = 0, $ which
verifies (D). Similarly, we can show that (E) is satisfied and thus
the functional central limit theorem holds.
\end{pf}

Let
\begin{eqnarray*}
G(t) & = & \alpha\pi_1 E P\bigl(\big|Z + \sqrt{n}\mu_1/\sigma_1\big|\ge t\bigr)- (1 - \alpha)(1 - \pi_1) P(|Z| \ge t) \\
& = & \alpha\pi_1 E P\bigl(\big|Z + \sqrt{n}|\mu_1|/\sigma_1\big|\ge t\bigr) - (1 -
\alpha)(1 - \pi_1) P(|Z| \ge t)
\end{eqnarray*}
and
%
%
\begin{equation}
t_1 = \inf\{t\dvtx  G(t) = 0\}.\ \label{t1}
\end{equation}
The following lemma is needed in the proof of consistency.
\begin{lemma} \label{l4.3}
Assume that $0< \pi_1 < 1 - \alpha$ and \eq{l4.2a} is satisfied. Then
%
%
\begin{equation}
G(t) \cases{
< 0 &\quad\mbox{for }$t < t_1$,\cr
=0 &\quad\mbox{for }$t = t_1$,\cr
>0 &\quad\mbox{for }$t > t_1$.
}
\label{l4.3a}
\end{equation}
Moreover, $G'(t_1) \ge \mathrm{e}^{-t_0^2/2} / \sqrt{2 \uppi}$.
\end{lemma}

\begin{pf}
We first observe that $0< t_1 \leq t_0$ by the fact
that $G(0) < 0,$ $G(t_0) > \mathrm{e}^{-t_0^2/2} > 0$ in \eq{l4.2-3} and
$G(t)$ is a continuous function.

To prove \eq{l4.3a}, it suffices to show that there exists a $t_2 >
t_1$ such that $G(t)$ is increasing in $[0, t_2]$ and decreasing in
$[t_2, \infty)$. To this end, consider the derivative of $G$:
%
%
\begin{eqnarray}\label{l4.3-1}
\hspace*{-50pt}G'(t) & = & -\alpha\pi_1E \bigl( \phi\bigl( t - \sqrt{n}|\mu_1|/\sigma_1\bigr)
+ \phi\bigl( t + \sqrt{n}|\mu_1|/ \sigma_1\bigr) \bigr) + 2 (1 - \alpha)(1 -\pi_1)\phi(t)\nonumber\\
& = & \frac{\mathrm{e}^{-t^2/2}}{\sqrt{2 \uppi}}
\biggl\{ - \alpha\pi_1
E \biggl(\exp\biggl(-\frac{n\mu_1^2}{2\sigma_1^2} +
\frac{\sqrt{n}|\mu_1|t}{\sigma_1} \biggr) +
\exp\biggl(-\frac{n\mu_1^2}{2\sigma_1^2} -
\frac{\sqrt{n}|\mu_1|t}{\sigma_1} \biggr) \biggr) \\
&&\hspace*{30pt}{}+ 2(1 - \alpha)(1 - \pi_1) \biggr\}.\nonumber
\end{eqnarray}
Let
\begin{eqnarray*}
H(t) &=& -\alpha\pi_1 E \biggl(\exp\biggl(-\frac{n\mu_1^2}{2\sigma_1^2} +
\frac{\sqrt{n}|\mu_1|t}{\sigma_1} \biggr)
\\
&&\hspace*{37pt}{}+\exp\biggl(-\frac{n\mu_1^2}{2\sigma_1^2} -
\frac{\sqrt{n}|\mu_1|t}{\sigma_1} \biggr) \biggr)
+ 2(1 - \alpha)(1 -
\pi_1).
\end{eqnarray*}
Then
%
%
\begin{eqnarray}\label{l4.3-2}
\hspace*{-15pt}H'(t) & = &- \alpha\pi_1E \biggl\{ \frac{\sqrt{n}|\mu_1| }{\sigma_1}
\exp\biggl( {\sqrt{n}|\mu_1| t \over\sigma_1} -
\frac{n\mu_1^2}{2\sigma_1^2} \biggr)  - \frac{\sqrt{n}|\mu_1| }{\sigma_1} \exp\biggl( -
\frac{\sqrt{n}|\mu_1| t}{\sigma_1} -
\frac{n\mu_1^2}{2\sigma_1^2} \biggr)
\biggr\} \nonumber\\[-8pt]\\[-8pt]
& = & - \alpha\pi_1 E \biggl\{
\frac{\sqrt{n}|\mu_1|}{\sigma_1}\mathrm{e}^{-\afrac{n\mu_1^2}{2\sigma_1^2}}
\biggl( \exp\biggl( \frac{\sqrt{n}|\mu_1| t}{\sigma_1} \biggr)
-\exp\biggl( -\frac{\sqrt{n}|\mu_1| t}{\sigma_1} \biggr) \biggr) \biggr\}
< 0\nonumber
\end{eqnarray}
for all $t>0$. Therefore, $H(t)$ is monotone decreasing. Taking into
account the facts that $H(0)> 0$ by assumption, $\pi_1 < 1 -
\alpha$ and $H(+\infty) < 0$, we conclude that $H(t)$ has only one
zero point, say, $t_2$. Moreover, $H(t) >0 $ for $t < t_2$ and
$H(t)< 0$ for $t>t_2$. This is also true for $G'(t)$, by \eq{l4.3-1}.
Hence, $G(t)$ is increasing for $t < t_2$ and decreasing for
$t>t_2$. Note that since $G(0) < 0, G(t_0) > 0$ and $G(+\infty) =
0$, we can see that $G(t)$ has a unique zero point $t_1$ and
$t_2>t_1$. Since $G(t)$ is increasing for $ 0< t < t_2$, we have
$G'(t_1)>0$. We now prove that $G'(t_1) \geq \mathrm{e}^{-t_0^2/2} /\sqrt{2 \uppi}$. It follows from the proof of \eq{l4.2-3} that
%
%
\begin{equation}
G(t_0) \geq \mathrm{e}^{-t_0^2/2}. \label{l4.3-10}
\end{equation}
Recalling that $G'(t)= {\mathrm{e}^{-t^2/2} \over\sqrt{2 \uppi}} H(t)$ and $H$
is decreasing, we have
%
\begin{eqnarray}\label{l4.3-11}
G(t_0) & = & G(t_0) - G(t_1) = \int_{t_1}^{t_0} G'(s )\,\mathrm{d}s
 \leq \int_{t_1}^{t_0} { \mathrm{e}^{-s^2/2} \over\sqrt{2 \uppi}} H(t_1)\,\mathrm{d}s\nonumber\\[-8pt]\\[-8pt]
& \leq& H(t_1) \bigl(1-\Phi(t_1)\bigr) \leq H(t_1) \mathrm{e}^{-t_1^2/2} = G'(t_1)\sqrt{ 2 \pi}.\nonumber
\end{eqnarray}
This proves $G'(t_1) \geq
\mathrm{e}^{-t_0^2/2}/\sqrt{2 \uppi}$.
\end{pf}

\subsection{\texorpdfstring{Proof of Theorem~\protect\ref{th3.1}}{Proof of Theorem 2.1}}

We now return to show our main theorem under dependence. Let
$\mathcal{H}= \{H_i, 1 \le i \le m\}$. To prove (i), following along
the same lines as the proof
of Lemma~\ref{th2.1}, we need to obtain the asymptotic distribution of
%
%
\begin{equation}
P(V \ge\alpha R) = P\Biggl(\sum_{i=1}^m \xi_i(t) \le0\Biggr), \label{4.5}
\end{equation}
where
\[
\xi_i(t) = \alpha I_{\{|T_i| \ge t\}} - I_{\{|T_i| \ge t\}}I_{\{H_i
= 0\}} = (\alpha+H_i - 1) I_{\{|T_i| \ge t\}}=[\alpha H_i - (1 -
\alpha
) (1 - H_i)]I_{\{|T_i| \ge t\}}.
\]
Note that
\[
P(|T_i|\geq t | \mathcal{H}) = (1-H_i) P(|T_i|\geq t|H_i=0) + H_i P(|T_i
|\geq t | H_i =1).
\]
Given $\mathcal{H}$, $\xi_i(t), 1 \leq i \leq m$, are independent random
variables. The conditional mean equals
\begin{eqnarray}\label{t3.1}
E\Biggl(\sum_{i=1}^m \xi_i | \mathcal{H}\Biggr)
& =& 
\sum_{i=1}^m \bigl\{ \alpha E\bigl( I_{\{H_i=0\}} | \mathcal{H}\bigr)
P(|T_i|\geq t | H_i=0)
+ \alpha E\bigl(I_{\{H_i=1\}} |\mathcal{H}\bigr) P(|T_i| \geq t |H_i=1)
\nonumber\\
&&\hspace*{16pt}{}-
E\bigl(I_{\{H_i=0\}} | \mathcal{H}\bigr) P( |T_i|\geq t | H_i=0) \bigr\}\nonumber\\
& = & \sum_{i=1}^m \{ \alpha(1-H_i) P(|T_i|\geq t |
H_i=0)
+ \alpha H_i P(|T_i| \geq t |H_i=1)\nonumber\\
&&\hspace*{16pt}{}- (1-H_i) P( |T_i|\geq t | H_i=0) \}\nonumber\\
& =& \alpha
\sum_{i=1}^m \{ H_i P(|T_i| \geq t
|H_i=1) \}- (1 - \alpha) 
\sum_{i=1}^m \{(1-H_i)
P(|T_i|\geq t |
H_1=0) \} \nonumber\\
& = & \alpha m_1 F_1(t) - (1 - \alpha) m_0 F_0(t). \nonumber
\end{eqnarray}

Next, we calculate the conditional variance of $\sum_{i=1}^m
\xi_i(t)$, given $\mathcal{H}$:
\begin{eqnarray*}
\operatorname{var}\Biggl(\sum_{i = 1}^m
\xi_i(t) | \mathcal{H}\Biggr)
& = & \operatorname{var}\Biggl(\sum_{i = 1}^m [\alpha H_i - (1- \alpha) (1- H_i)]
I_{\{|T_i| \ge t | \mathcal{H}\}}\Biggr)\\
& = & \sum_{i = 1}^{m} \bigl(\alpha^2 H_i + (1 - \alpha)^2 (1 -
H_i)\bigr)\operatorname{var}\bigl(I_{\{|T_i| \ge t | \mathcal{H}\}}\bigr) \\
& = & \alpha^2 m_1 F_1(t)\bigl(1 - F_1(t)\bigr) + (1 - \alpha)^2 m_0 F_0(t)\bigl(1- F_0(t)\bigr).
\end{eqnarray*}
From (\ref{mu}) and (\ref{sigma}),
\[
\frac{\mu_m(t)}{\sigma_m(t)} =
\sqrt{m}\frac{\mu_m(t)/m}{\sqrt{\sigma^2_m(t)/m}}.
\]
By the fact that $m_1/m \rightarrow\pi_1$ a.s., we have
%
%
\begin{equation}
\mu_m(t)/m \rightarrow\alpha\pi_1 F_1(t) - (1 - \alpha) (1 - \pi_1)
F_0(t)\qquad \mbox{a.s.} \label{muxi}
\end{equation}
and
%
%
\begin{eqnarray}
\sigma^2_m(t)/m &\rightarrow&\alpha^2 \pi_1 F_1(t)\bigl(1 - F_1(t)\bigr)\nonumber
\\[-8pt]\\[-8pt]
&&{}+ (1 -\alpha)^2 (1 - \pi_1)F_0(t)\bigl(1 - F_0(t)\bigr) = q^2(t) \qquad\mbox{a.s.},\nonumber
\end{eqnarray}
which is smaller than $\operatorname{var}(\xi_1(t))$, due to the fact that
\[
\operatorname{var} X = E(\operatorname{var}(X|Y)) + \operatorname{var}(E(X|Y))
\]
for any two random variables $X$ and $Y$. By (\ref{l4.2-1}), we can
see that the critical value defined at (\ref{tnm}) is bounded. Thus,
conditional on $\mathcal{H}$, we can use the functional central limit
theorem on
$\sum_{i=1}^m \xi_i(t)/{\sqrt{m}}$, by virtue of Lemma~\ref{lem7.8}.
The limit is a Gaussian process with continuous sample paths. Hence,
\begin{eqnarray*}
P\Biggl(\sum_{i = 1}^{m}\xi_i(t)\le0\Biggr) 
& = & E \bigl(E 1_{\{\sum_{i = 1}^{m}\xi_i(t)/\sqrt{m} \le0\}} \big|
\mathcal{H}\bigr) \\
& = & E \Biggl\{P \Biggl(\sum_{i = 1}^{m} \xi_i/\sqrt{m} - \sum_{i = 1}^{m}
E (\xi_i |\mathcal{H})/\sqrt{m} \le
\frac{-\sum_{i = 1}^{m} E (\xi_i |\mathcal{H})\sigma_m(t)}{\sqrt
{m}\sigma_m(t)}
\bigg| \mathcal{H}\Biggr) \Biggr\} \\
& \le& E \Biggl\{P \Biggl(\sum_{i = 1}^{m} \xi_i/\sqrt{m} - \sum_{i = 1}^{m}
E (\xi_i |\mathcal{H})/\sqrt{m} \le
\frac{-\sum_{i = 1}^{m} E (\xi_i |\mathcal{H})}{\sigma_m(t)}\frac
{\sigma
_m(t)}{\sqrt{m}} \bigg| \mathcal{H}\Biggr) \Biggr\} \\
& \le& E \bigl\{P \bigl(N(0, 1)q(t) \le
-z_{\gamma}q(t) \bigr) \bigr\}\\
& \rightarrow& P\bigl(N(0, 1) \le-z_{\gamma}\bigr) = \gamma\qquad \mbox{as }m
\rightarrow\infty.
\end{eqnarray*}
This proves (\ref{tnm}).

(ii) can be proven similarly. The characteristic function method can
be used to prove (iii).

\subsection{\texorpdfstring{Proof of Theorem~\protect\ref{th3.2}}{Proof of Theorem 2.2}}

We first prove (i), and (ii) follows along the same lines as the
independent case, plus a conditional argument.
Without loss of generality, we use $T_1$ as a representative that comes from
the alternative. We have to show that
%
%
\begin{equation}
|\hat{t}_{n,m} - t_{n,m}| = \mathrm{o}(1)\qquad \mbox{a.s.} \label{l4.4a}
\end{equation}

We first prove
that
%
%
\begin{equation}
|\hat{t}_{n,m} - t_1| = \mathrm{o}(1)\qquad \mbox{a.s.}, \label{l4.4-3}
\end{equation}
where $t_1$ is defined as in (\ref{t1}). It suffices to show that
for any $\varepsilon>0$,
%
%
\begin{equation}
{\sqrt{m} \nu_m(t_1 + \varepsilon) \over\tau_m(t_1+\varepsilon)}
\geq z_\gamma
\label{l4.4-1}
\end{equation}
and
%
%
\begin{equation}
{\sqrt{m}\nu_m(s) \over\tau_m(s)} < z_\gamma\qquad  \mbox{for all }
 s \le t_1 - \varepsilon. \label{l4.4-2}
\end{equation}

Recall that $\hat{p}_m(t) = \frac{1}{m}\sum_{i = }^{m}I_{\{|T_i|\ge t\}}$.
Given $\cH$, by the uniform law of the iterated logarithm
(see, e.g., \cite{r13}),
\begin{eqnarray*}\label{dudley}
\hat{p}_m(t) - \frac{1}{m}\sum_{i = 1}^{m}\{(1 - H_i)F_0(t) +
H_i F_1(t) \} = \mathrm{o}(m^{-1/2}(\log\log m)^{1/2})\qquad\mbox{a.s.}
\end{eqnarray*}
By the strong law of large number,
\begin{equation}
\frac{1}{m}\sum_{i = 1}^{m}\{(1 - H_i)F_0(t) + H_i F_1(t) \}
\rightarrow (1 - \pi_1)F_0(t) + \pi_1 F_1(t) \qquad  \mbox{a.s.}
\end{equation}
So
\begin{eqnarray*}
\hat{p}_m(t) \rightarrow (1 - \pi_1)F_0(t) + \pi_1 F_1(t) \qquad \mbox{a.s.}
\end{eqnarray*}
Recall that
\[
\nu_m(t) = \alpha \hat{p}_m(t) - 2 (1 - \hat{\pi}_1)\bar{\Phi}(t).
\]
By \eq{l4.5a}, our strong consistent estimate $\hat{\pi}_1$
described in Section~\ref{sec23} and the continuous mapping theorem, we have
\begin{equation}\label{l4.4-4}
\sup_t \big|\nu_m(t) - \bigl\{ \alpha \bigl((1 - \pi_1)F_0(t) + \alpha \pi_1
F_1(t)\bigr) - (1 - \pi_1)P(|Z| \ge t)\bigr\}\big| \to 0 \qquad \mbox{a.s.},
\end{equation}
which, together with \eq{4.03} and the definition of $G$, implies
that
\begin{equation}\label{l4.4-5}
\sup_{0 \leq t \leq 1+ t_0} |\nu_m(t) - G(t)| \to 0  \qquad\mbox{a.s.}
\end{equation}
%
In particular, since $G(t_1+\varepsilon) >0$ for $0< \varepsilon<
t_2-t_1$, we
have
%
%
\begin{equation}
\nu_m(t_1+\varepsilon) \geq G(t_1+\varepsilon)/2\qquad \mbox{a.s.} \label{l4.4-6}
\end{equation}
for sufficiently large $m$ and, therefore, $\sqrt{m}\nu_m(t_1 +
\epsilon) \ge z_{\gamma} \tau_m(t_1 + \epsilon)$. This proves
\eq{l4.4-1}.

Similarly, since $G(t)$ is increasing and $G(t_1-\varepsilon) < 0$, we have
%
%
\begin{equation}
\max_{s \leq t_1 -\varepsilon} \nu_m(s) \leq G(t_1-\varepsilon)/2\qquad \mbox{a.s.}
\label{l4.4-7}
\end{equation}
for sufficiently large $m$. Hence, \eq{l4.4-2} holds. This proves
\eq{l4.4-3}.

Following the same lines as the proof of \eq{l4.4-3}, we have
%
%
\begin{equation}
|t_{n,m} - t_1| = \mathrm{o}(1). \label{l4.4-8}
\end{equation}
This completes the proof of \eq{l4.4a}.

For $k$-FWER, let $\eta_0$ be the number that satisfies $P(\operatorname{Poiss}(\eta_0) \ge k) \le\gamma$. Let $t_{0,m}=t_{n,m}^{k\mbox{-}\mathit{FWER}}$ and $t_m=\hat{t}_{m,n}^{k\mbox{-}\mathit{FWER}}$.
Thus, by definition, $t_{0,m}$ is the $t$ that satisfies $(1 - \pi_1) m
F_{\mathrm{o}}(t) = \eta_0$
and $t_{m}$ is the $t$ that satisfies $2(1 - \hat{\pi}_1) m \bar
{\Phi
}(t) = \eta_0.$
We then have
$\frac{(1 - \pi_1)F_0(t_{0,m})}{(1 - \hat{\pi}_1)2 \bar{\Phi
}(t_{m})} =
1$, which implies that
\begin{eqnarray*}
&&\frac{F_0(t_{0,m})}{2\bar{\Phi}(t_m)}= \frac{1 -\hat{\pi}_1}{1 - \pi_1}  =  1 + \mathrm{o}_P(1)\\
&&\quad \Longrightarrow\quad\frac{\bar{\Phi}(t_{0,m})}{\bar{\Phi}(t_m)}\bigl(1 +\mathrm{O}(n^{-1/2})\bigr) = 1 +\mathrm{o}_P(1) \\
&&\quad\Longrightarrow\quad\frac{\bar{\Phi}(t_{0,m})}{\bar{\Phi}(t_m)} =  1 + \mathrm{o}_P(1)\\
&&\quad\Longrightarrow\quad\frac{t_{m}}{t_{0,m}} \mathrm{e}^{-t_{0,m}^2/2 + t_m^2/2} =  1 + \mathrm{o}_P(1)\\
&&\quad\Longrightarrow\quad R \mathrm{e}^{-t_{0,m}^2/2 + R^2 t_{0,m}^2/2} = R \mathrm{e}^{-(1 - R^2)t_{0,m}^2/2}=1 + \mathrm{o}_P(1).
\end{eqnarray*}
Hence, $R = t_m/t_{0,m} \to1$ in probability.
Thus,
\begin{eqnarray*}
t_{0,m}^2 - t_{m}^2& = & \mathrm{o}_P(1) \quad\Longrightarrow\quad
|t_{0,m} - t_m| = \frac{\mathrm{o}_P(1)}{1 + |t_{0,m} + t_{m}|} = \mathrm{O}_p({(\log m)}^{-1/2})
\end{eqnarray*}
since $t_m = \mathrm{o}_P(n^{1/6})$ and $\log m = \mathrm{o}(n^{1/3}).$

\subsection{\texorpdfstring{Proof of Theorem~\protect\ref{t2.2}}{Proof of Theorem 2.4}}

In this section, we give the proof of the rate of convergence for the
i.i.d.~case by using the one-sample $t$-statistic. Let $p(t) = P(|T_1|
\geq t)$ and let
\[
\hat{p}_m(t) = { 1 \over m} \sum_{i=1}^m I_{\{|T_i|\geq t\}}.
\]
By the Glivenko--Cantelli theorem,
%
%
\begin{equation}
\sup_{t} |\hat{p}_m(t) - p(t)| \to0 \qquad\mbox{a.s.} \label{4.20}
\end{equation}
and, by the Donsker theorem,
%
%
\begin{equation}
\sup_{t} |\hat{p}_m(t) - p(t)| = \mathrm{O}(m^{-1/2}) \qquad\mbox{in
probability}. \label{4.21}
\end{equation}
By the uniform law of the iterated logarithm,
%
%
\begin{equation}
\sup_{t} |\hat{p}_m(t) - p(t)| = \mathrm{O}(m^{-1/2}(\log\log m)^{1/2})
\qquad\mbox{a.s.} \label{4.21-2}
\end{equation}
We define strong consistent estimators of $E\xi_1(t)$ and
$\operatorname{var}(\xi_1(t))$ by $\nu_m(t)$ and $\tau_m^2(t)$,
respectively, where
%
%
\begin{equation}
\nu_m(t) = \alpha\hat{p}_m(t) - (1 - \pi_1)P(|Z| \ge t) \label{nu}
\end{equation}
and
%
%
\begin{eqnarray}\label{tau-1}
\tau_m^2(t)& = & \alpha^2\hat{p}_m(t)\bigl(1 - \hat{p}_m(t)\bigr)
+ 2 \alpha(1- \pi_1)\hat{p}_m(t)P(|Z| \ge t)\nonumber\\[-8pt]\\[-8pt]
&&{}+ (1 -\pi_1)P(|Z| \ge t)\bigl(1 - 2\alpha- (1 - \pi_1)P(|Z| \ge t)\bigr).\nonumber
\end{eqnarray}
We now define an estimator of $t_{n,m}$ by
%
%
\begin{equation}
\hat{t}_{n,m} = \inf\biggl\{t\dvtx  {\frac{\sqrt{m}\nu_m(t)}{\tau_m(t)}} \ge
z_{\gamma}\biggr\}. \label{4.22}
\end{equation}

For \mbox{FDTP}, we have to show that
%
%
\begin{equation}
|\hat{t}_{n,m} - t_{n,m}| = \mathrm{O}\biggl(\frac{1}{\sqrt{n}} + \biggl(\frac{\log\log
m}{m}\biggr)^{1/2}\biggr)\qquad \mbox{a.s.} \label{l4.4b}
\end{equation}
and
%
%
\begin{equation}
|\hat{t}_{n,m} - t_{n,m}| = \mathrm{O}(n^{-1/2} + m^{-1/2})
\qquad \mbox{in probability}.
\label{l4.4c}
\end{equation}

Below, we prove \eq{l4.4b} and \eq{l4.4c}. We will show that
%
%
\begin{eqnarray}
|\hat{t}_{n,m} - t_1| = \mathrm{O}\biggl(\biggl(\frac{1}{n}\biggr)^{1/2} + \biggl(\frac{\log\log
m}{m}\biggr)^{1/2}\biggr) \qquad\mbox{a.s.}, \label{l4.5-11}
\\\label{l4.5-12}
|t_{n,m} - t_1| = O\biggl(\biggl(\frac{1}{n}\biggr)^{1/2} + \biggl(\frac{\log\log
m}{m}\biggr)^{1/2}\biggr) \qquad\mbox{a.s.}
\end{eqnarray}
By the uniform law of the iterated logarithm,
%
%
\begin{equation}
\sup_t |\hat{p}_m(t) - p(t)| = \mathrm{O}\biggl(\biggl(\frac{\log\log m}{m}\biggr)^{1/2}\biggr)
 \qquad\mbox{a.s.} \label{lil}
\end{equation}
Therefore, we have
%
%
\begin{equation}
\sup_t \big|v_m(t) - [\alpha p(t) - (1 - \pi_1)P(|Z| \ge t)]\big| =
\mathrm{O}\biggl(\biggl(\frac{\log\log m}{m}\biggr)^{1/2}\biggr) \qquad\mbox{a.s.} \label{l4.5-13}
\end{equation}
Note that
\begin{eqnarray*}
&& \alpha p(t) - (1 - \pi_1)P(|Z| \ge t) - G(t)\\
&&\quad =
\alpha(1 -\pi_1)\bigl( P(|T_1| \ge t |H_1 = 0) - P(|Z| \ge t)\bigr) \\
&&\qquad{} +  \alpha\pi_1 \bigl(P(|T_1|
\ge t |H_1 = 1) -E P\bigl(\big|Z + \sqrt{n}\mu_1/\sigma_1\big| \ge t\bigr) \bigr).\\
\end{eqnarray*}
From \eq{l4.5a}, we obtain
%
%
\begin{equation}
P(|T_1| \ge t|H_1 = 0) - P(|Z| \ge t) = \mathrm{O}\biggl(\frac{1}{\sqrt{n}}\biggr)
\qquad\mbox{a.s.} \label{null}
\end{equation}
and %
%
%
\begin{equation}
P(|T_1| \ge t|H_1 = 1) - E P\bigl(\big|Z + \sqrt{n}\mu_1/\sigma_1\big| \ge t\bigr)
=\mathrm{O}\biggl(\frac{1}{\sqrt{n}}\biggr) \qquad\mbox{a.s.} \label{alt}
\end{equation}
Thus,
we have
%
%
\begin{equation}
\sup_t\big|\alpha p(t) - (1 - \pi_1)P(|Z| \ge t) - G(t)\big| =
\mathrm{O}\biggl(\frac{1}{\sqrt{n}}\biggr) \qquad\mbox{a.s.}
\end{equation}
Taking into account \eq{l4.5-13}, we have
%
%
\begin{equation}
\sup_t|v_m(t) - G(t)| \leq c_2 \biggl(\frac{1}{\sqrt{n}} + \biggl(\frac{\log
\log m}{m}\biggr)^{1/2}\biggr) \qquad\mbox{a.s.} \label{l4.5-14}
\end{equation}
for some constant $0< c_2 < \infty$. Below, we show that there exists
a finite constant $c_3>0$ such that
%
%
\begin{equation}
t_1 - c_3\biggl(\frac{1}{\sqrt{n}} + \biggl(\frac{\log\log m}{m}\biggr)^{1/2}\biggr) <
\hat{t}_{n,m} < t_1 + c_3\biggl(\frac{1}{\sqrt{n}} + \biggl(\frac{\log\log
m}{m}\biggr)^{1/2}\biggr). \label{l4.5-15}
\end{equation}

Recalling \eq{l4.5-14}, we have, for $\epsilon= c_3(\frac{1}{\sqrt
{n}} +
(\frac{\log\log m}{m})^{1/2})$, that
\begin{eqnarray*}
v_m(t_1 + \epsilon) & \ge& G(t_1 + \epsilon) -
c_2\biggl(\frac{1}{\sqrt{n}}
+ \biggl(\frac{\log\log m}{m}\biggr)^{1/2}\biggr)\\
& = & G(t_1) + \epsilon G'(t_1 + \theta_1) - c_2\biggl(\frac{1}{\sqrt{n}}
+
\biggl(\frac{\log\log m}{m}\biggr)^{1/2}\biggr) \\
& \ge& c_1 \epsilon- c_2\biggl(\frac{1}{\sqrt{n}} + \biggl(\frac{\log\log
m}{m}\biggr)^{1/2}\biggr) > 2 \biggl(\frac{\log\log m}{m} \biggr)^{1/2},
\end{eqnarray*}
provided that $c_3$ is chosen large enough: here, $0 \le\theta_1 \le
\epsilon$ and we have used Lemma~\ref{l4.3}. For sufficiently large $m$,
we have
\[
\sqrt{m}v_m(t_1 + \epsilon) > \tau_m(t_1 + \epsilon) z_\gamma.
\]
This proves
that
\[
\hat{t}_{n,m} - t_1 \leq c_3 \biggl(\biggl(\frac{1}{n}\biggr)^{1/2} + \biggl(\frac{\log\log
m}{m}\biggr)^{1/2}\biggr) \qquad\mbox{a.s.}
\]
Similarly, we have
\[
\hat{t}_{n,m} - t_1 \geq- c_3 \biggl(\biggl(\frac{1}{n}\biggr)^{1/2} +
\biggl(\frac{\log\log m}{m}\biggr)^{1/2}\biggr) \qquad\mbox{a.s.}
\]

This proves \eq{l4.5-11}.

Following the same line of proof, we have
\[
| t_{n,m} - t_1| = \mathrm{O}\biggl(\frac{1}{\sqrt{n}} + \biggl(\frac{\log\log
m}{m}\biggr)^{1/2}\biggr) \qquad\mbox{a.s.}
\]
If we use
%
%
\begin{equation}
\sup_t |\hat{p}_m(t) - p(t)| = \mathrm{O}(m^{-1/2}) \qquad\mbox{in
probability,}\ \label{lil-p}
\end{equation}
based on the Donsker theorem instead of \eq{lil}, using the same line
of the proof of the a.s.~convergence rate, we can obtain the rate of
convergence in
probability, which is
\[
|\hat{t}_{n,m} - t_{n,m}| = \mathrm{O}(n^{-1/2} + m^{-1/2})
\qquad\mbox{in probability}.
\]
This completes the proof of \eq{l4.4b}.

Similarly, the critical value for FDR control is bounded, due to the
fact that
\[
E P\biggl(\bigg|Z + \frac{\sqrt{n}\mu_1}{\sigma_1}\bigg| \ge t\biggr) \le1.
\]
By (\ref{4.21}), (\ref{4.21-2}), (\ref{null}) and (\ref{alt}), we have
\begin{eqnarray*}
\sup_t\bigg|\frac{m_0F_0(t)}{m_0F_0(t) + m_1F-1(t)} - \frac{2(1 - \pi_1)\bar{\Phi}(t)}{\hat{p}_m(t)}\bigg|
&=& \mathrm{O}\biggl(n^{-1/2} +\biggl(\frac{\log\log m}{m}\biggr)^{1/2}\biggr) \qquad\mbox{a.s.},\\
\sup_t\bigg|\frac{m_0F_0(t)}{m_0F_0(t) + m_1F-1(t)} - \frac{2(1 - \pi_1)\bar{\Phi}(t)}{\hat{p}_m(t)}\bigg|
&=& \mathrm{O}\bigl(n^{-1/2} +(m)^{-1/2}\bigr)\qquad \mbox{in probability.}
\end{eqnarray*}
Noting that $2(1 - \pi_1)\bar{\Phi}(t)/[2(1 - \pi_1)\bar{\Phi}(t) +
EP(|Z + \sqrt{n}\mu_1/\sigma_1|\ge t)]$ is a monotone decreasing
continuous function with respect to $t$, combined with the definitions
of ($t_{n,m}^{\mathit{fdr}}$) and ($\hat{t}_{n,m}^{\mathit{fdr}}$), (\ref{rate-1}) and
(\ref{rate-2}) hold.

The proof of $k$-FWER is the same as that given in Theorem~\ref{th3.2}.

\subsection{\texorpdfstring{Proof of Theorem~\protect\ref{th2.3}}{Proof of Theorem 3.1}}

For the two-sample $t$-statistic, the only part we
need to show is the boundedness of $t_{n,m}$ under independence, which
will imply the boundedness in the general dependence case, as happens
with the one-sample $t$-statistic. The remaining results follows along
the same lines as the proof in the one sample $t$-statistic setting.
Based on Lemma~\ref{l5.6}
below, plus \eq{th2.3}, and using the same line of proof as in the
one-sample $t$-statistic case, the boundedness of $t_{n,m}$ holds for
two-sample $t$-statistics.

The proof of the boundedness of $t_{n,m}$ is based on the following
asymptotic distribution of $T_i^*$ under the alternative hypothesis.
\begin{lemma} \label{l5.6}
Suppose that $X, X_1, \ldots, X_{n_1}$ are independent and identically
distributed random variables from a population with mean $\mu_1$ and
variance $\sigma_1^2$, and $Y, Y_1, \ldots, Y_{n_2}$ are independent and
identically distributed random variables from another population
with mean $\mu_2$ and variance $\sigma_2^2.$ Assume the sampling
processes are independent of each other. Also, assume that there are
$0 < c_1 \le c_2 < \infty$ such that $c_1 \le n_1/n_2 \le c_2$. Let
%
%
\begin{equation}
T^*= \frac{\bar{X} - \bar{Y} 
}{\sqrt{s_1^2/n_1 + s_2^2/n_2}},
\end{equation}
where
%
%
\begin{eqnarray}
\bar{X} &=& \frac{1}{n_1}\sum_{i = 1}^{n_1} X_i,\qquad
\bar{Y} = \frac{1}{n_2}\sum_{i = 1}^{n_2} Y_i,
\\
s_1^2 &=& \frac{1}{n_1 - 1} \sum_{i = 1}^{n_1}(X_i - \bar{X})^2
\quad\mbox{and}\quad
s_2^2 = \frac{1}{n_2 - 1} \sum_{i = 1}^{n_2}(Y_i -\bar{Y})^2.
\end{eqnarray}
If $EX^4< \infty$ and $EY^4 < \infty$, then
%
%
\begin{equation}
P(|T^*| \ge t) 
= P\biggl(\bigg|Z + \frac{\mu_1 - \mu_2}{\sqrt{\sigma_1^2/n_1 +
\sigma_2^2/n_2}}\bigg| \ge t\biggr) \bigl(1 + \mathrm{o}(1)\bigr),
\end{equation}
uniformly in $t = \mathrm{o}(n^{1/6})$, where $n = \max{\{n_1, n_2\}}.$
\end{lemma}

\begin{pf}
The proof of this lemma is very similar to the proof of Lemma~\ref{l4.1} and so we omit the details.
\end{pf}

\subsection{\texorpdfstring{Proof of Theorem~\protect\ref{th2.4}}{Proof of Theorem 3.2}}

This follows the same arguments as in the one-sample $t$-statistic case,
by virtue of Lemma~\ref{l5.6}.

\subsection{\texorpdfstring{Proof of Theorem~\protect\ref{in-2}}{Proof of Theorem 3.3}}

When we plug in an estimator of $P(|T_i^*| \ge t)$,
\[
\hat{p}_m(t) = \frac{1}{m}\sum_{i = 1}^{m}I_{\{|T_i^*| \ge t\}},
\]
the proof of the two-sample $t$-statistic case follows along the same
lines as
its one-sample counterpart, except that we have to show the rate of
convergence under the alternative hypothesis for the two-sample
$t$-statistic. This follows from the following lemma, which completes
the proof of Theorem~\ref{in-2}.
\begin{lemma}\label{l5.7}
Let $X, X_1, \ldots, X_{n_1}$ be i.i.d.~random variables from a
population with mean $\mu_1$ and variance $\sigma_1^2$, and $Y, Y_1,
\ldots, Y_{n_2}$ be i.i.d.~random variables from another population
with mean $\mu_2$ and variance $\sigma_2^2.$ The sampling processes
are assumed to be independent of each other. Assume that there are
$0 < c_1 \le c_2 < \infty$ such that $c_1 \le n_1/n_2 \le c_2.$ Let
$T^*$ be defined as in Lemma~\textup{\ref{l5.6}}.
If $E|X|^4 < \infty$ and $E|Y|^4 < \infty$, then
%
%
\begin{eqnarray}\label{5.70}
&&\bigg|P(T^* \le x) - \Phi\biggl(x - \frac{\mu_1 - \mu_2}{\sqrt{\sigma_1^2/n_1+\sigma_2^2/n_2}}\biggr)\bigg|\nonumber
\\[-8pt]\\[-8pt]
&&\quad\le\frac{K(1 + |x|)}{(1 + | x- \vfrac{\mu_1 -\mu_2}{\sqrt{\sigma_1^2/n_1 + \sigma_2^2/n_2}}|)\sqrt{\min\{n_1,n_2\}}},\nonumber
\end{eqnarray}
where $K$ is a finite constant that may depend on $\sigma_1^2,
\sigma_2^2, E|X|^3, E|Y|^3, EX^4$ and $EY^4$.
\end{lemma}

\begin{pf}
Without loss of generality, we assume that $n_1 = b_1 n$, $n_2
=b_2 n$, $b_1 + b_2 = 1$ with $b_1 >0$ and $b_2 >0$. Note that
\begin{eqnarray*}
P(T^* \le x) & = & P\biggl(\frac{\bar{X} - \mu_1 - (\bar{Y} -\mu_2)}{\sqrt{s_1^2/n_1 + s_2^2/n_2}}
+ \frac{\mu_1 -\mu_2}{\sqrt{s_1^2/n_1 + s_2^2/n_2}} \le x\biggr)\\
& = & P\biggl(\frac{\bar{X} - \mu_1 - (\bar{Y} -
\mu_2)}{\sqrt{\sigma_1^2/n_1 + \sigma_2^2/n_2}} + \frac{\mu_1 -
\mu_2}{\sqrt{\sigma_1^2/n_1 + \sigma_2^2/n_2}} \le
x\frac{\sqrt{s_1^2/n_1 + s_2^2/n_2}}{\sqrt{\sigma_1^2/n_1 +
\sigma_2^2/n_2}}\biggr)\\
& \le& P\biggl(\frac{\bar{X} - \mu_1 - (\bar{Y} -
\mu_2)}{\sqrt{\sigma_1^2/n_1 + \sigma_2^2/n_2}} \le x - \frac{\mu_1
- \mu_2}{\sqrt{\sigma_1^2/n_1 + \sigma_2^2/n_2}} + x\bigg|\frac{s_1^2/n_1
+ s_2^2/n_2}{\sigma_1^2/n_1 + \sigma_2^2/n_2} - 1\bigg|\biggr),
\end{eqnarray*}
where we\vspace*{-6pt} make use of \eq{00}. We now apply \eq{l4.05a} with $\xi_i =
\frac{(X_i - \mu_1)/n_1}{\sqrt{\sigma_1^2/n_1 + \sigma_2^2/n_2}}$
for $1 \le i \le n_1$ and $\xi_i = -\frac{(Y_i -
\mu_2)/n_2}{\sqrt{\sigma_1^2/n_1 + \sigma_2^2/n_2}}$ for $n_1 + 1
\le i \le n_1 + n_2$. Let
\begin{eqnarray*}
z &=& x - \frac{\mu_1 - \mu_2}{\sqrt{\sigma_1^2/n_1 +\sigma_2^2/n_2}},\qquad
\\
\Delta&=& -x\bigg|\frac{s_1^2/n_1 +s_2^2/n_2}{\sigma_1^2/n_1 + \sigma_2^2/n_2} - 1\bigg|,
\\
\Delta_i &=& -x\bigg|\frac{s_{1,i}^2/n_1 +
s_2^2/n_2}{\sigma_1^2/n_1 + \sigma_2^2/n_2} - 1\bigg|
\end{eqnarray*}
for $1 \le i \le n_1$, and
\[
\Delta_i = -x\bigg|\frac{s_{1}^2/n_1 + s_{2,i}^2/n_2}{\sigma_1^2/n_1 +
\sigma_2^2/n_2} - 1\bigg|
\]
for $n_1 + 1 \le i \le n_1 + n_2$, where $s^2_{1,i}$ is defined as
$s^2_1$ with 0 replacing
$X_i$ and $s^2_{2,i}$ is defined as $s^2_2$ with 0 replacing $Y_i$.
Noting that
\begin{eqnarray*}
\frac{s_1^2/n_1 + s_2^2/n_2}{\sigma_1^2/n_1 + \sigma_2^2/n_2} - 1 &
= & \frac{1}{\sigma_1^2/n_1 + \sigma_2^2/n_2}[(s_1^2 -
\sigma_1^2)/n_1 + (s_2^2 - \sigma_2^2)/n_2],
\end{eqnarray*}
%
we have, by \eq{l4.5-5}, that
\[
E\bigg|\frac{s_1^2/n_1 + s_2^2/n_2}{\sigma_1^2/n_1 + \sigma_2^2/n_2} -1\bigg|^2 \le K \frac{EX^4 + EY^4}{n}.
\]
For $1 \le i \le n_1$,
\begin{eqnarray*}
&& E\biggl(\frac{s_1^2/n_1 + s_2^2/n_2}{\sigma_1^2/n_1 +\sigma_2^2/n_2}- \frac{s_{1i}^2/n_1 + s_2^2/n_2}{\sigma_1^2/n_1 +\sigma_2^2/n_2}\biggr)^2 \\
&&\quad =  \frac{1}{n_1^2(\sigma_1^2/n_1 + \sigma_2^2/n_2)^2}E(s_1^2 - s_{1i}^2)^2\le\frac{K EX^4}{n^2},
\end{eqnarray*}
by \eq{l4.5-6}. Similarly, for $n_1 + 1 \le i \le n_1 + n_2$, we have
\begin{eqnarray*}
&& E\biggl(\frac{s_1^2/n_1 + s_2^2/n_2}{\sigma_1^2/n_1 +
\sigma_2^2/n_2}
- \frac{s_{1}^2/n_1 + s_{2i}^2/n_2}{\sigma_1^2/n_1 +
\sigma_2^2/n_2}\biggr)^2
 =  \frac{1}{n_2^2(\sigma_1^2/n_1 +
\sigma_2^2/n_2)^2}E(s_2^2 - s_{2i}^2)
\le\frac{K EY^4}{n^2}.
\end{eqnarray*}
It follows that
\begin{eqnarray*}
\|\Delta\|_2 & \le& K
\frac{|x|\sqrt{EX^4 + EY^4}}{\sqrt{n}},\\
P\biggl(|\Delta|> \frac{|z|+1}{3}\biggr) & \le& K\frac{E|\Delta|}{|z| + 1} \le
K \frac{\|\Delta\|_2}{|z|+1} \le
K \frac{|x|\sqrt{EX^4 + EY^4}}{\sqrt{n}(|z|+1)},\\
\sum_{i = 1}^{n}(E\xi_i^2)^{1/2}\bigl(E(\Delta-\Delta_i)^2\bigr)^{1/2} & \le
& K \frac{\sqrt{(\sigma_1^2 + \sigma_2)(EX^4 + EY^4)}}{\sqrt{n}},\\
\sum_{i=1}^{n}E|\xi_i|^3 & \le& K \frac{E|X|^3 + E|Y|^3}{\sqrt{n}}.
\end{eqnarray*}
%
Therefore, by \eq{l4.05a},
\begin{eqnarray*}
&&\bigg|P\biggl(\frac{\bar{X} - \mu_1 - (\bar{Y} - \mu_2)}{\sqrt{\sigma
_1^2/n_1 +
\sigma_2^2/n_2}}  \le x - \frac{\mu_1 -
\mu_2}{\sqrt{\sigma_1^2/n_1 + \sigma_2^2/n_2}} + x\bigg|\frac{s_1^2/n_1
+ s_2^2/n_2}{\sigma_1^2/n_1 + \sigma_2^2/n_2} - 1\bigg|\biggr)\\
&&\quad{} -  \Phi\biggl(x - \frac{\mu_1 - \mu_2}{\sqrt{\sigma_1^2/n_1 +
\sigma_2^2/n_2}}\biggr)\bigg| \le K \frac{1 + |x|}{(1 + | x - \vfrac{\mu_1 -
\mu_2}{\sqrt{\sigma_1^2/n_1 + \sigma_2^2/n_2}}|)\sqrt{n}}.
\hspace{8cm}
\end{eqnarray*}
Similarly,
\begin{eqnarray*}
P(T^* \le x) &= & P\biggl(\frac{\bar{X} - \mu_1 - (\bar{Y} -
\mu_2)}{\sqrt{s_1^2/n_1 + s_2^2/n_2}} + \frac{\mu_1 -
\mu_2}{\sqrt{s_1^2/n_1 + s_2^2/n_2}} \le
x\biggr) \\
& \ge& P\biggl(\frac{\bar{X} - \mu_1 - (\bar{Y} -
\mu_2)}{\sqrt{\sigma_1^2/n_1 + \sigma_2^2/n_2}} \le x - \frac{\mu_1
- \mu_2}{\sqrt{\sigma_1^2/n_1 + \sigma_2^2/n_2}} - x\bigg|\frac{s_1^2/n_1
+ s_2^2/n_2}{\sigma_1^2/n_1 + \sigma_2^2/n_2} - 1\bigg|\biggr)
\end{eqnarray*}
and
\begin{eqnarray*}
&&\bigg|P\biggl(\frac{\bar{X} - \mu_1 - (\bar{Y} - \mu_2)}{\sqrt
{\sigma_1^2/n_1
+ \sigma_2^2/n_2}}  \le x - \frac{\mu_1 -
\mu_2}{\sqrt{\sigma_1^2/n_1 + \sigma_2^2/n_2}} - x\bigg|\frac{s_1^2/n_1 +
s_2^2/n_2}{\sigma_1^2/n_1 + \sigma_2^2/n_2} - 1\bigg|\biggr) \\
&&\quad{} -  \Phi\biggl(x -
\frac{\mu_1 - \mu_2}{\sqrt{\sigma_1^2/n_1 + \sigma_2^2/n_2}}\biggr)\bigg|
\le
K \frac{1 + |x|}{(1 + | x - \vfrac{\mu_1 -
\mu_2}{\sqrt{\sigma_1^2/n_1 + \sigma_2^2/n_2}}|)\sqrt{n}}.
\hspace{8cm}
\end{eqnarray*}

This proves \eq{5.70}. \end{pf}

\end{appendix}

\section*{Acknowledgements}
The authors are grateful to the Associate Editor and two anonymous
reviewers for valuable comments and suggestions which improved the
paper. The second author was partially supported by grants CA075142 and
CA142538 from the U.S. National Institutes of Health and grant
DMS-0904184 from the U.S. National Science Foundation.


\printhistory

\end{document}